\documentclass[a4paper,11pt]{article}

\usepackage{latexsym}
\usepackage{amssymb}
\usepackage{amsmath}
\usepackage{color}
\usepackage{graphicx}
\usepackage{natbib}
\usepackage{amsfonts}
\usepackage{amscd}
\usepackage{vmargin}
\usepackage[T1]{fontenc}
\usepackage[english]{babel}
\usepackage[latin1]{inputenc}

\newcommand{\p}{\mathbb{P}}
\newcommand{\C}{\mathbb{C}}
\newcommand{\R}{\mathbb{R}}
\newcommand{\Z}{\mathbb{Z}}

\newcommand{\N}{\mathbb{N}}
\newcommand{\E}{\mathbb{E}}

\newcommand{\M}{\mathcal{M}}
\newcommand{\F}{\mathcal{F}}
\newcommand{\pen}{\textrm{pen}}
\newcommand{\Cov}{\textrm{Cov}}
\newcommand{\Lip}{\textrm{Lip}}
\newcommand{\Var}{\textrm{Var}}

\newtheorem{theo}{Theorem}[section]

\newtheorem{prop}[theo]{Proposition}
\newtheorem{lemma}[theo]{Lemma}

\newenvironment{proof}[2][Proof]{ 
\begin{trivlist} 
\item[\hskip \labelsep {\bfseries \itshape \textit{#1  #2}}]~\\ \indent}
{\end{trivlist}}

\begin{document}

%\chapter[Model selection for mixing data]{Adaptive density estimation of stationary $\beta$-mixing and $\tau$-mixing processes}
%\begin{flushright}
%{\small{{\it Mon modèle, c'est moi-même! Je suis mon meilleur modèle parce que je connais mes erreurs, mes qualités, mes victoires et mes défaites. Si je passe mon temps à prendre un autre modèle comme modèle, comment veux-tu que ce modèle puisse modeler dans la bonne ligne?\\
%JC Vandamme}}}
%\end{flushright}
\title{Adaptive density estimation of stationary $\beta$-mixing and $\tau$-mixing processes.}
\date{}
\author{Matthieu Lerasle\footnote{Institut de Math\'ematiques (UMR 5219), INSA de Toulouse, 
Universit\'e de Toulouse, France}}
\maketitle

\hspace{1cm}\begin{minipage}{12cm}
\begin{center}
 Abstract:
\end{center}
{\small We propose an algorithm to estimate the common density $s$ of a stationary process $X_1,...,X_n$. We suppose that the process is either $\beta$ or $\tau$-mixing. We provide a model selection procedure based on a generalization of Mallows' $C_p$ and we prove oracle inequalities for the selected estimator under a few prior assumptions on the collection of models and on the mixing coefficients. We prove that our estimator is adaptive over a class of Besov spaces, namely, we prove that it achieves the same rates of convergence as in the i.i.d framework.}
\end{minipage}
\vspace{0.5cm}



\section{Introduction}
We consider the problem of estimating the unknown density $s$ of $P$, the law of a random variable $X$, based on the observation of $n$ (possibly) dependent data $X_1,...,X_n$ with common law $P$. We assume that $X$ is real valued, that $s$ belongs to $L^2(\mu)$ where $\mu$ denotes the Lebesgue measure on $\R$ and that $s$ is compactly supported, say in $[0,1]$. Throughout the chapter, we consider least-squares estimators $\hat{s}_m$ of $s$ on a collection $(S_m)_{m\in\M_n}$ of linear subspaces of $L^2(\mu)$. Our final estimator is chosen through a model selection algorithm.\\
Model selection has received much interest in the last decades. When its final goal is prediction, it can be seen more generally as the question of choosing between the outcomes of several prediction algorithms. With such a general formulation, a very natural answer is the following. First, estimate the prediction error for each model, that is $\|s-\hat{s}_m\|_2^2$. Then, select the model which minimizes this estimate.\\
It is natural to think of the empirical risk as an estimator of the prediction error. This can fail dramatically, because it uses the same data for building predictors and for comparing them, making these estimates strongly biased for models involving a number of parameters growing with the sample size.\\
In order to correct this drawback, penalization's methods state that a good choice can be made by minimizing the sum of the empirical risk (how do algorithms fit the data) and some complexity measure of the algorithms (called the penalty). This method was first developped in the work of Akaike \cite{Ak70} and \cite{Ak73} and Mallows \cite{Ma73}.\\
In the context of density estimation, with independent data, Birg\'e $\&$ Massart \cite{BM97} used penalties of order $L_nD_m/n$, where $D_m$ denotes the dimension of $S_m$ and $L_n$ is a constant depending on the complexity of the collection $\M_n$. They used Talagrand's inequality (see for example Talagrand \cite{Ta96} for an overview) to prove that this penalization procedure is efficient {\it{i.e.}} the integrated quadratic risk of the selected estimator is asymptotically equivalent to the risk of the oracle (see Section 2 for a precise definition). They also proved that the selected estimator achieves adaptive rates of convergence over a large class of Besov spaces. Moreover, they showed that some methods of adaptive density estimation like the unbiased cross validation (Rudemo \cite{Ru82}) or the hard thresholded estimator of Donoho {\it{et al.}} \cite{DJKP} can be viewed as special instances of penalized projection estimators. \\
More recently, Arlot \cite{Ar08} introduced new measures of the quality of penalized least-squares estimators (PLSE). He proved pathwise oracle inequalities, that is deviation bounds for the PLSE that are harder to prove but more informative from a practical point of view (see also Section 2 for details).\\
When the process $(X_i)_{i=1,...,n}$ is $\beta$-mixing (Rozanov $\&$ Volkonskii \cite{RV59} and Section 2), Talagrand's inequality can not be used directly. Baraud {\it {et al.}} \cite{BCV01} used Berbee's coupling lemma (see Berbee (\cite{Be79}) and Viennet's covariance inequality (Viennet \cite{Vi97}) to overcome this problem and build model selection procedure in the regression problem. Then Comte $\&$ Merlev\`ede \cite{CM02} used this algorithm to investigate the problem of density estimation for a $\beta$-mixing process. They proved that under reasonable assumptions on the collection $\M_n$ and on the coefficients $\beta$, one can recover the results of Birg\'e $\&$ Massart \cite{BM97} in the i.i.d. framework.\\
The main drawback of those results is that many processes, even simple Markov chains are not $\beta$-mixing. For instance, if $(\epsilon_i)_{i\geq 1}$ is iid with marginal $\mathcal{B}(1/2)$, then the stationary solution $(X_i)_{i\geq 0}$ of the equation 
\begin{equation}\label{andrews}
X_n=\frac{1}{2}(X_{n-1}+\epsilon_n),\;X_0\;\textrm{independent of}\;(\epsilon_i)_{i\geq 1}
\end{equation}
is not $\beta$-mixing (Andrews \cite{An84}). More recently, Dedecker $\&$ Prieur \cite{DP05} introduced new mixing-coefficients, in particular the coefficients $\tau$, $\tilde{\phi}$ and $\tilde{\beta}$ and proved that many processes like (\ref{andrews}) happen to be $\tau$, $\tilde{\phi}$ and $\tilde{\beta}$-mixing. They proved a coupling lemma for the coefficient $\tau$ and covariance inequalities for $\tilde{\phi}$ and $\tilde{\beta}$. Gannaz $\&$ Wintenberger \cite{GW09} used the covariance inequality to extend the result of Donoho {\it{et al.}} \cite{DJKP} for the wavelet thresholded estimator to the case of $\tilde{\phi}$-mixing processes. They recovered (up to a $\log(n)$ factor) the adaptive rates of convergence over Besov spaces. \\
In this article, we first investigate the case of $\beta$-mixing processes. We prove a pathwise oracle inequality for the PLSE. We extend the result of Comte $\&$ Merlev\`ede \cite{CM02} under weaker assumptions on the mixing coefficients. Then, we consider $\tau$-mixing processes. The problem is that the coupling result is weaker for the coefficient $\tau$ than for $\beta$. Moreover, in order to control the empirical process we use a covariance inequality that is harder to handle. Hence, the generalization of the procedure of Baraud {\it{et al.}} \cite{BCV01} to the framework of $\tau$-mixing processes is not straightforward. We recover the optimal adaptive rates of convergence over Besov spaces (that is the same as in the independent framework) for $\tau$-mixing processes, which is new as far as we know.\\
The chapter is organized as follows. In Section 2, we give the basic material that we will use throughout the chapter. We recall the definition of some mixing coefficients and we state their properties. We define the penalized least-squares estimator (PLSE). Sections 3 and 4 are devoted to the statement of the main results, respectively in the $\beta$-mixing case and in the $\tau$-mixing case. In Section 5, we derive the adaptive properties of the PLSE. Finally, Section 6 is devoted to the proofs. Some additional material has been reported in the Appendix in Section 7.
\section{Preliminaries}
\subsection{Notation.}
Let $(\Omega,\mathcal{A},\p)$ be a probability space. Let $\mu$ be the Lebesgue measure on $\R$, let $\left\|.\right\|_p$ be the usual norm on $L^p(\mu)$ for $1\leq p\leq \infty$. For all $y\in\R^l$, let $\left|y\right|_l=\sum_{i=1}^l|y_i|$. Denote by $\lambda_{\kappa}$ the set of $\kappa$-Lipschitz functions, {\it{i.e.}} the functions $t$ from $(\R^l,\left|.\right|_l)$ to $\R$ such that $\Lip(t)\leq \kappa$ where
$$\Lip(t)=\sup\left\{\frac{|t(x)-t(y)|}{|x-y|_l},x,y\in \R^l,x\neq y\right\}\leq \kappa.$$ 
Let $BV$ and $BV_1$ be the set of functions $t$ supported on $\R$ satisfying respectively $\left\|t\right\|_{BV}<\infty$ and $\left\|t\right\|_{BV}\leq 1$ where 
$$\left\| t\right\|_{BV}=\sup_{n\in\N^*}\sup_{-\infty<a_1<...<a_n<\infty}|t(a_{i+1})-t(a_i)|.$$
\subsection{Some measures of dependence.} 
\subsubsection{Definitions and assumptions}
Let $Y=(Y_1,...,Y_l)$ be a random variable defined on $(\Omega,\mathcal{A},\p)$ with values in $(\R^l,\left|.\right|_l)$. Let $\M$ be a $\sigma$-algebra of $\mathcal{A}$. Let $\p_{Y|\M}$, $\p_{Y_1|\M}$ be conditional distributions of $Y$ and $Y_1$ given $\M$, let $\p_Y$, $\p_{Y_1}$ be the distribution of $Y$ and $Y_1$ and let $F_{Y_1|\M}$, $F_{Y_1}$ be distribution functions of $\p_{Y_1|\M}$ and $P_{Y_1}$. Let $\mathcal{B}$ be the Borel $\sigma$-algebra on $(\R^l,\left|.\right|_l)$. Define now 
\begin{eqnarray*}
\beta(\M,\sigma(Y))&=&\E\left(\sup_{A\in\mathcal{B}}|\p_{Y|\M}(A)-\p_Y(A)|\right),\\
\tilde{\beta}(\M,Y_1)&=&\E\left(\sup_{x\in\R}\left|F_{Y_1|\M}(x)-F_{Y_1}(x)\right|\right),\\
\textrm{and if}\; \E(|Y|)<\infty,\; \tau(\M,Y)&=&\E\left(\sup_{t\in\lambda_1}|\p_{Y|\M}(t)-\p_Y(t)|\right).
\end{eqnarray*}
The coefficient $\beta(\M,\sigma(Y))$ is the mixing coefficient introduced by Rozanov $\&$ Volkonskii \cite{RV59}. The coefficients $\tilde{\beta}(\M,Y_1)$ and $\tau(\M,Y)$ have been introduced by Dedecker $\&$ Prieur \cite{DP05}.\\
Let $(X_k)_{k\in\Z}$ be a stationary sequence of real valued random variables defined on $(\Omega,\mathcal{A},\p)$. For all $k\in\N^*$, the coefficients $\beta_k$, $\tilde{\beta}_k$ and $\tau_k$ are defined by 
\begin{equation*}
\beta_k=\beta(\sigma(X_i,  i\leq 0),\sigma(X_i,  i\geq k)),\;\tilde{\beta}_k=\sup_{j\geq k}\lbrace\tilde{\beta}(\sigma(X_p,  p\leq 0),X_j)\rbrace.
\end{equation*}
If $\E(|X_1|)<\infty$, for all $k\in\N^*$ and all $r\in\N^*$, let 
\begin{equation*}
\tau_{k,r}=\max_{1\leq l\leq r}\frac{1}{l}\sup_{k\leq i_1<..< i_l}\lbrace\tau(\sigma(X_p,  p\leq 0),(X_{i_1},...,X_{i_l}))\rbrace,\;\tau_k=\sup_{r\in\N^*}\tau_{k,r}.
\end{equation*} 
Moreover, we set $\beta_0=1$. In the sequel, the processes of interest are either $\beta$-mixing or $\tau$-mixing, meaning that, for $\gamma=\beta$ or $\tau$, the $\gamma$-mixing coefficients $\gamma_k\rightarrow 0$ as $k\rightarrow +\infty$. For $p\in \left\{1,2\right\}$, we define $\kappa_p$ as: 
\begin{equation}\label{1}
\kappa_p=p\sum_{l=0}^{\infty}l^{p-1}\beta_l,
\end{equation}
where $0^0=1$, when the series are convergent. Besides, we consider two kinds of rates of convergence to $0$ of the mixing coefficients, that is for $\gamma=\beta$ or $\tau$,\\
${\bf{[AR]}}$ arithmetical $\gamma$-mixing with rate $\theta$ if there exists some $\theta>0$ such that $\gamma_k\leq (1+k)^{-(1+\theta)}$ for all $k$ in $\N$,\\
${\bf{[GEO]}}$ geometrical $\gamma$-mixing with rate $\theta$ if there exists some $\theta>0$ such that $\gamma_k\leq e^{-\theta k}$ for all $k$ in $\N$.
\subsubsection{Properties}
{\bf{Coupling}}\\
Let $X$ be an $\R^l$-valued random variable defined on $(\Omega,\mathcal{A},\p)$ and let $\M$ be a $\sigma$-algebra. Assume that there exists a random variable $U$ uniformly distributed on $[0,1]$ and independent of $\M\vee\sigma(X)$. There exist two $\M\vee\sigma(X)\vee\sigma(U)$-measurable random variables $X_1^*$ and $X_2^*$ distributed as $X$ and independent of $\M$ such that
\begin{equation}\label{coup1}
\beta(\M,\sigma(X))=\p(X\neq X_1^*)\;\textrm{and}\;
\end{equation}
\begin{equation}\label{coup2}
\tau(\M,X)=\E\left(|X-X_2^*|_l\right).
\end{equation}
Equality (\ref{coup1}) has been established by Berbee \cite{Be79}, Equality (\ref{coup2}) has been established in Dedecker $\&$ Prieur \cite{DP05}, Section 7.1.\\
{\bf{Covariance inequalities}}\\
Let $X,Y$ be two real valued random variables and let $f,h$ be two measurable functions from $\R$ to $\C$. Then, there exist two measurable functions $b_1:\R\rightarrow \R$ and $b_2:\R\rightarrow \R$ with $\E \left(b_1(X)\right)=\E( b_2(Y))=\beta(\sigma(X),\sigma(Y))$ such that, for any conjugate $p,q\geq1$ (see Viennet \cite{Vi97} Lemma 4.1)
\begin{equation*}\label{covB}
|\Cov(f(X),h(Y))|\leq 2\E^{1/p} \left( |f(X)|^pb_1(X)\right)\E^{1/q}(|h(Y)|^{q}b_2(Y)).
\end{equation*}
There exists a random variable $b(\sigma(X),Y)$ such that $\E (b(\sigma(X),Y))=\tilde{\beta}(\sigma(X),Y)$ and such that, for all Lipschitz functions $f$ and all $h$ in $BV$ (Dedecker $\&$ Prieur \cite{DP05} Proposition 1)
\begin{equation}\label{covT}
|\Cov(f(X),h(Y))|\leq \left\| h\right\|_{BV}\E\left( |f(X)|b(\sigma(X),Y)\right)\leq \left\| h\right\|_{BV}\left\|f\right\|_{\infty}\tilde{\beta}(\sigma(X),Y).
\end{equation}
{\bf{Comparison results}}\\
Let $(X_k)_{k\in\Z}$ be a sequence of identically distributed real random variables. If the marginal distribution satisfies a concentration's condition $|F_X(x)-F_X(y)|\leq K|x-y|^a$ with $a\leq 1$, $K>0$, then (Dedecker {\it{et al.}} \cite{DDLLLP} Remark 5.1 p 104)
$$\tilde{\beta}_k\leq 2K^{1/(1+a)}\tau_{k,1}^{a/(a+1)}\leq 2K^{1/(1+a)}\tau_{k}^{a/(a+1)}.$$
In particular, if $\p_X$ has a density $s$ with respect to the Lebesgue measure $\mu$ and if $s\in L^2(\mu)$, we have from Cauchy-Schwarz inequality 
$$|F_X(x)-F_X(y)|=|\int{\bf{1}}_{[x,y]}sd\mu|\leq \left\|s\right\|_2\left(\int{\bf{1}}_{[x,y]}d\mu\right)^{1/2}=\left\|s\right\|_2|x-y|^{1/2},$$
thus 
\begin{equation*}\label{comp}
\tilde{\beta}_k\leq 2\left\|s\right\|^{2/3}_2\tau_k^{1/3}.
\end{equation*}
In particular, for any arithmetically ${\bf{[AR]}}$ $\tau$-mixing process with rate $\theta>2$, we have \begin{equation}\label{cvg}
\tilde{\beta}_k\leq2\left\|s\right\|^{2/3}_2(1+k)^{-(1+\theta)/3}.
\end{equation}
\subsubsection{Examples}
Examples of $\beta$-mixing and $\tau$-mixing sequences are well known, we refer to the books of Doukhan \cite{Do94} and Bradley \cite{Br07} for examples of $\beta$-mixing processes and to the book of Dedecker {\it et. al} \cite{DDLLLP} or the articles of Dedecker $\&$ Prieur \cite{DP05}, Prieur \cite{Pr07}, and Comte {\it et. al} \cite{CDT08} for examples of $\tau$-mixing sequences. One of the most important example is the following: a stationary, irreducible, aperiodic and positively recurent Markov chain $(X_i)_{i\geq 1}$ is $\beta$-mixing. However, many simple Markov chains are not $\beta$-mixing but are $\tau$-mixing. For instance, it is known for a long time that if $(\epsilon_i)_{i\geq 1}$ are i.i.d Bernoulli $\mathcal{B}(1/2)$, then a stationary solution $(X_i)_{i\geq 0}$ of the equation
$$X_n=\frac{1}{2}(X_{n-1}+\epsilon_n),\;X_0\;\textrm{independent of}\;(\epsilon_i)_{i\geq 1}$$ 
is not $\beta$-mixing since $\beta_k=1$ for any $k\geq 1$ whereas $\tau_k\leq 2^{-k}$ (see Dedecker $\&$ Prieur \cite{DP05} Section 4.1). Another advantage of the coefficient $\tau$ is that it is easy to compute in many situations (see Dedecker $\&$ Prieur \cite{DP05} Section 4). 

\subsection{Collections of models}
We observe $n$ identically distributed real valued random variables $X_1,...,X_n$ with common density $s$ with respect to the Lebesgue measure $\mu$. We assume that $s$ belongs to the Hilbert space $L^2(\mu)$ endowed with norm $\left\|.\right\|_2$. We consider an orthonormal system $\left\lbrace \psi_{j,k}\right\rbrace_{(j,k)\in \Lambda}$ of $L_2(\mu)$ and a collection of models $(S_m)_{m\in\M_n}$ indexed by subsets $m\subset\Lambda$ for which we assume that the following assumptions are fulfilled:\\
{\bf {[$M_1$]}} for all $m\in\M_n$, $S_m$ is the linear span of $\left\lbrace \psi_{j,k}\right\rbrace_{(j,k)\in m}$ with finite dimension $D_m=|m|\geq 2$ and $N_n=\max_{m\in\M_n} D_m$ satisfies $N_n\leq n$;\\
{\bf {[$M_2$]}} there exists a constant $\Phi$ such that 
$$\forall m,m'\in\M_n,\forall t\in S_m,\forall t'\in S_{m'}, \|t+t'\|_{\infty}\leq \Phi\sqrt{\textrm{dim}(S_m+S_{m'})}\|t+t'\|_2;$$
{\bf {[$M_3$]}} $D_m\leq D_{m'}$ implies that $m\subset {m'}$ and so $S_m\subset S_{m'}$.\\
As a consequence of Cauchy-Schwarz inequality, we have
\begin{equation}\label{norme}
\left\|\sum_{(j,k)\in m\cup m'}\psi_{j,k}^2\right\|_{\infty}=\sup_{t\in S_m+S_{m'}, t\neq 0}\frac{\|t\|_{\infty}^2}{\|t\|_2^2}
\end{equation}
see Birg\'e $\&$ Massart \cite{BM97} p 58. Three examples are usually developed as fulfilling this set of assumptions:\\
{\bf{[T]}} trigonometric spaces: $\psi_{0,0}(x)=1$ and for all $j\in \N^*$, $\psi_{j,1}(x)=\cos(2\pi jx)$, $\psi_{j,2}(x)=\sin(2\pi jx)$. $m=\{(0,0),(j,1),(j',2),\;1\leq j,j' \leq J_m\}$ and $D_m=2J_m+1$;\\
{\bf{[P]}} regular piecewise polynomial spaces: $S_m$ is generated by $r$ polynomials $\psi_{j,k}$ of degree $k=0,...,r-1$ on each subinterval $[(j-1)/J_m,j/J_m]$ for $j=1,...,J_m$, $D_m=rJ_m$, $\M_n=\left\{m=\{(j,k),\;j=1,...,J_m,\;k=0,...,r-1\},1\leq J_m\leq [n/r]\right\}$;\\
{\bf{[W]}} spaces generated by dyadic wavelet with regularity $r$ as described in Section 4.\\
For a precise description of those spaces and their properties, we refer to Birg\'e $\&$ Massart \cite{BM97}.

\subsection{The estimator}
Let $(X_n)_{n\in\Z}$ be a real valued stationary process and let $P$ denote the law of $X_0$. Assume that $P$ has a density $s$ with respect to the Lebesgue measure $\mu$ and that $s\in L_2(\mu)$. Let $(S_m)_{m\in\M_n}$ be a collection of models satisfying assumptions {\bf {[$M_1$]}}-{\bf {[$M_3$]}}. We define $S_n=\cup_{m\in\M_n}S_m$, $s_m$ and $s_n$ the orthogonal projections of $s$ onto $S_m$ and $S_n$ respectively, let $\p$ be the joint distribution of the observations $(X_n)_{n\in\Z}$ and let $\E$ be the corresponding expectation. We define the operators $P_n$, $P$ and $\nu_n$ on $L^2(\mu)$ by
$$P_nt=\frac{1}{n}\sum_{i=1}^nt(X_i),\;Pt=\int t(x)s(x)d\mu(x),\;\nu_n(t)=(P_n-P)t.$$
All the real numbers that we shall introduce and which are not indexed by $m$ or $n$ are fixed constants. In order to define the penalized least-squares estimator, let us consider on $\R\times S_n$ the contrast function $\gamma(x,t)=-2t(x)+\left\|t\right\|_2^2$ and its empirical version $\gamma_n(t)=P_n\gamma(.,t)$. Minimizing $\gamma_n(t)$ over $S_m$ leads to the classical projection estimator $\hat{s}_m$ on $S_m$. Let $\hat{s}_n$ be the projection estimator on $S_n$. Since $\left\lbrace \psi_{j,k}\right\rbrace_{(j,k)\in m}$ is an orthonormal basis of $S_m$ one gets 
$$\hat{s}_m=\sum_{(j,k)\in m}(P_n\psi_{j,k})\psi_{j,k}\;\textrm{and}\;\gamma_n(\hat{s}_m)=-\sum_{(j,k)\in m}(P_n\psi_{j,k})^2.$$
Now, given a penalty function $\pen : \M_n\rightarrow \R^+$, we define a selected model $\hat{m}$ as any element
\begin{equation}\label{selmodel}
\hat{m}\in\textrm{arg}\min_{m\in\M_n}\left(\gamma_n(\hat{s}_m)+\pen(m)\right)
\end{equation}
and a PLSE is defined as any $\tilde{s}\in S_{\hat{m}}\subset S_n$ such that 
\begin{equation}\label{PLSE}
\gamma_n(\tilde{s})+\pen(\hat{m})=\inf_{m\in\M_n}\left(\gamma_n(\hat{s}_m)+\pen(m)\right).
\end{equation}

\subsection{Oracle inequalities}
An ideal procedure for estimation chooses an oracle 
$$m_o\in\textrm{Arg}\min_{m\in\M_n}\lbrace\Vert s-\hat{s}_m\Vert_2\rbrace.$$
An oracle depends on the unknown $s$ and on the data so that it is unknown in practice. In order to validate our procedure, we try to prove:\\
-non asymptotic oracle inequalities for the PLSE: 
\begin{equation}\label{oracleC}
\E\left(\left\| s-\tilde{s}\right\|_2^2\right)\leq L\inf_{m\in\M_n}\lbrace\E\left(\Vert s-\hat{s}_m\Vert_2^2+R(m,n)\right)\rbrace,
\end{equation}
for some constant $L\geq 1$ (as close to $1$ as possible) and a remainder term $R(m,n)\geq 0$ possibly random, and small compared to $\E\left(\left\| s-\tilde{s}\right\|_2^2\right)$ if possible. This inequality compares the risk of the PLSE with the best deterministic choice of $m$. Since $\hat{m}$ is random, we prefer to prove a stronger form of oracle inequality :
\begin{equation}\label{oracleR}
\E\left(\left\| s-\tilde{s}\right\|_2^2\right)\leq L\E\left(\inf_{m\in\M_n}\lbrace\left\| s-\hat{s}_{m}\right\|_2^2+R(m,n)\rbrace\right),
\end{equation}
or, when it is possible, deviation bounds for the PLSE: 
\begin{equation}\label{oracleP}
\p\left(\left\| s-\tilde{s}\right\|_2^2> L\inf_{m\in\M_n}\left(\left\| s-\hat{s}_{m}\right\|_2^2+R(m,n)\right)\right)\leq c_n,
\end{equation}
where typically $c_n\leq C/n^{1+\gamma}$ for some $\gamma>0$. Inequality (\ref{oracleP}) proves that, asymptotically, the risk $\left\| s-\tilde{s}\right\|_2^2$ is almost surely the one of the oracle. Let 
$$\Omega=\left\{\left\| s-\tilde{s}\right\|_2^2> L\inf_{m\in\M_n}\left(\left\| s-\hat{s}_{m}\right\|_2^2+R(m,n)\right)\right\}.$$
We have
$$\E\left(\left\| s-\tilde{s}\right\|_2^2\right)= \E\left(\left\| s-\tilde{s}\right\|_2^21_{\Omega}\right)+\E\left(\left\| s-\tilde{s}\right\|_2^21_{\Omega^c}\right).$$
It is clear that $\E\left(\left\| s-\tilde{s}\right\|_2^21_{\Omega^c}\right)\leq L\E\left(\inf_{m\in\M_n}\left\{\left\| s-\hat{s}_{m}\right\|_2^2+R(m,n)\right\}\right).$ Moreover, we have $\|s-\tilde{s}\|^2=\|s-s_{\hat{m}}\|^2+\|s_{\hat{m}}-\tilde{s}\|^2\leq \|s\|^2+\Phi^2D_{\hat{m}}\leq \|s\|^2+\Phi^2n$, thus, when (\ref{oracleP}) holds, we have 
$$ \E\left(\left\| s-\tilde{s}\right\|_2^21_{\Omega^c}\right)\leq (\|s\|^2+\Phi^2n)c_n\leq \frac{C}{n^{\gamma}}.$$
Therefore, inequality (\ref{oracleP}) implies
$$\E\left(\left\|s-\tilde{s}\right\|_2^2\right)\leq \E\left(\inf_{m\in\M_n}\lbrace\left\| s-\hat{s}_{m}\right\|_2^2+R(m,n)\rbrace\right)+\frac{C}{n^{\gamma}}.$$
We can derive from these inequalities adaptive rates of convergence of the PLSE on Besov spaces (see Birg\'e $\&$ Massart \cite{BM97} for example). In order to achieve this goal, we only have to prove a weaker form of oracle inequality where the remainder term $R(m,n)\leq LD_m/n$ for some constant $L$, for all the models $m$ with sufficiently large dimension. This will be detailed in Section 5.

\section{Results for $\beta$-mixing processes}
From now on, the letters $\kappa$, $L$ and $K$, with various sub- or supscripts, will denote some constants which may vary from line to line. One shall use $L_{.}$ to indicate more precisely the dependence on various quantities, especially those which are related to the unknown $s$.\\
In this section, we give the following theorem for $\beta$-mixing sequences. It can be seen as a pathwise version of Theorem 3.1 in Comte $\&$ Merlev\`ede \cite{CM02}.
\begin{theo}\label{betathm}
Consider a collection of models satisfying {\bf {[$M_1$]}}, {\bf {[$M_2$]}} and {\bf {[$M_3$]}}. Assume that the process $(X_n)_{n\in \Z}$ is strictly stationary  and arithmetically {\bf {[AR]}} $\beta$-mixing with mixing rate $\theta>2$ and that its marginal distribution admits a density $s$ with respect to the Lebesgue measure $\mu$, with $s\in L_2(\mu)$. \\
Let $\kappa_1$ be the constant defined in (\ref{1}) and let $\tilde{s}$ be the PLSE defined by (\ref{PLSE}) with 
$$\pen(m)=\frac{K\Phi^2\kappa_1D_m}{n},\;\textrm{where}\; K>4.$$
Then, for all $\kappa>2$ there exist $c_0>0,L_s>0$, $\gamma_1>0$ and a sequence $\epsilon_n\rightarrow 0$,  such that
\begin{equation}\label{Beta}
\p\left(\left\|\tilde{s}-s\right\|_2^2>(1+\epsilon_n)\inf_{m\in\M_n,D_m\geq c_0(\log n)^{\gamma_1}}\left( \left\|s-s_{m}\right\|_2^2+\pen(m) \right)\right)\leq L_s\frac{(\log n)^{(\theta+2)\kappa}}{n^{\theta/2}}.
\end{equation}
\end{theo}
{\bf Remark:} 
The term $K\Phi^2\kappa_1$ is the same as in Theorem 3.1 of Comte $\&$ Merlev\`ede \cite{CM02} but with a constant $K>4$ instead of $320$. The main drawback of this result is that the penalty term involves the constant $\kappa_1$ which is unknown in practice. However, Theorem \ref{betathm} ensures that penalties proportional to the linear dimension of $S_m$ lead to efficient model selection procedures. Thus we can use this information to apply the slope heuristic algorithm introduced by Birg\'e $\&$ Massart \cite{BM07} in a Gaussian regression context and generalized by Arlot $\&$ Massart \cite{AM08} to more general M-estimation frameworks. This algorithm calibrates the constant in front of the penalty term when the shape of an ideal penalty is available. The result of Arlot $\&$ Massart is proven for independent sequences, in a regression framework, but it can be generalized to the density estimation framework, for independent as well as for $\beta$ or $\tau$ dependent data. This result is beyond the scope of this chapter and will be proved in chapter 4. \\
We have to consider the infimum in equation (\ref{Beta}) over the models with sufficiently large dimensions. However, as noted by Arlot \cite{Ar08} (Remark 9 p 43), we can take the infimum over all the models in (\ref{Beta}) if we add an extra term in (\ref{Beta}). More precisely, we can prove that, with probability larger than $1-L_s(\log n)^{(\theta+2)\kappa}/n^{\theta/2}$
\begin{equation}
\left\|\tilde{s}-s\right\|_2^2\leq (1+\epsilon_n)\inf_{m\in\M_n}\left( \left\|s-\hat{s}_{m}\right\|_2^2+\pen(m)\right)+L\frac{(\log n)^{\gamma_2}}{n},
\end{equation}
where $L>0$ and $\gamma_2>0$.\\
{\bf Remark :}
The main improvement of Theorem \ref{betathm} is that it gives an oracle inequality in probability, with a deviation bound of order $o(1/n)$ as soon as $\theta>2$ instead of $\theta >3$ in Comte $\&$ Merlev\`ede \cite{CM02}. Moreover, we do not require $s$ to be bounded to prove our result.\\
{\bf Remark:}
When the data are independent, the proof of Theorem \ref{betathm} can be used to obtain that the estimator $\tilde{s}$ chosen with a penalty term of order $K\Phi D_m/n$ satisfy an oracle inequality as (\ref{Beta}). The main difference would be that $\kappa_1=1$, thus it can be used without a slope heuristic (even if this algorithm can be used also in this context to optimize the constant K) and the control of the probability would be $L_se^{-\ln(n)^2/C_s}$ for some constants $L_s,C_s$ instead of  $L_s(\log n)^{(\theta+2)}\kappa n^{-\theta/2}$ in our theorem.

\section{Results for $\tau$-mixing sequences}
In order to deal with $\tau$-mixing sequences, we need to specify the basis $(\psi_{j,k})_{(j,k)\in\Lambda}$.
\subsection{Wavelet basis}
Throughout this section, $r$ is a real number, $r\geq 1$ and we work with an $r$-regular orthonormal multiresolution analysis of $L_2(\mu)$, associated with a compactly supported scaling function $\phi$ and a compactly supported mother wavelet $\psi$. Without loss of generality, we suppose that the support of the functions $\phi$ and $\psi$ is an interval $[A_1,A_2)$ where $A_1$ and $A_2$ are integers such that $A_2-A_1=A\geq 1$. Let us recall that $\phi$ and $\psi$ generate an orthonormal basis by dilatations and translations. \\
For all $k\in \Z$ and $j\in \N^*$, let $\psi_{0,k}:\;x\rightarrow \sqrt{2}\phi(2x-k)$ and $\psi_{j,k}:\; x\rightarrow 2^{j/2}\psi(2^jx-k)$. The family $\left\{(\psi_{j,k})_{j\geq 0,k\in\Z}\right\}$ is an orthonormal basis of $L_2(\mu)$. Let us recall the following inequalities: for all $p\geq 1$, let $K_p=(\sqrt{2}\|\phi\|_p)\vee\|\psi\|_p$, $K_L=(2\sqrt{2}\Lip(\phi))\vee\Lip(\psi)$, $K_{BV}=AK_L$.\\
Then for all $j\geq 0$, we have $\|\psi_{j,k}\|_{\infty}\leq K_{\infty}2^{j/2},$
\begin{eqnarray}
\left\|\sum_{k\in\Z}|\psi_{j,k}|\right\|_{\infty}&\leq &AK_{\infty}2^{j/2}\label{inf}\\
\Lip(\psi_{j,k})&\leq &K_L2^{3j/2},\label{lip}\\
\|\psi_{j,k}\|_{BV}&\leq &K_{BV}2^{j/2}.\label{bv}
\end{eqnarray}
We assume that our collection $(S_m)_{m\in\M_n}$ satisfies the following assumption:\\
{\bf{[W]}} dyadic wavelet generated spaces: let $J_n=[\log(n/2(A+1))/\log(2)]$ and for all $J_m=1,...,J_n$, let 
$$m=\{(0,k),-A_2<k<2-A_1\}\cup\{(j,k),\; 1\leq j\leq J_m, \; -A_2< k< -A_1+2^j\}$$ 
and $S_m$ the linear span of $\{\psi_{j,k}\}_{(j,k)\in m}$. In particular, we have $D_m=(A-1)(J_m+1)+2^{J_m+1}$ and thus $2^{J_m+1}\leq D_m\leq (A-1)(J_m+1)+2^{J_m+1}\leq A2^{J_m+1}$.
\subsection{The $\tau$-mixing case}
The following result proves that we keep the same rate of convergence for the PLSE based on $\tau$-mixing processes.
\begin{theo}\label{tauthm}
Consider the collection of models {\bf{[W]}}. Assume that $(X_n)_{n\in \Z}$ is strictly stationary and arithmetically {\bf {[AR]}} $\tau$-mixing with mixing rate $\theta>5$ and that its marginal distribution admits a density $s$ with respect to the Lebesgue measure $\mu$. Let $\tilde{s}$  be the PLSE defined by (\ref{PLSE}) with 
$$\pen(m)=KAK_{\infty}K_{BV}\left(\sum_{l=0}^{\infty}\tilde{\beta}_l\right)\frac{D_m}{n}, \; \textrm{where}\; K\geq8.$$
Then there exist constants $c_0>0,\gamma_1>0$ and a sequence $\epsilon_n\rightarrow 0$ such that
\begin{equation}\label{tau}
\E\left(\|\tilde{s}-s\|_2^2\right)\leq (1+\epsilon_n)\left(\inf_{m\in \M_n,\;D_m \geq c_0(\log n)^{\gamma_1}}\|s-s_m\|_2^2+\pen(m)\right).
\end{equation}
\end{theo}
{\bf Remark :}
As in Theorem \ref{betathm}, the penalty term involves an unknown constant and we have a condition on the dimension of the models in (\ref{tau}). However, the slope heuristic can also be used in this context to calibrate the constant and a careful look at the proof shows that we can take the infimum over all models $m\in \M_n$ provided that we increase the constant $K$ in front of the penalty term. Our result allows to derive rates of convergence in Besov spaces for the PLSE that correspond to the rates in the i.i.d. framework (see Proposition \ref{minimax}).\\
{\bf Remark :}
Theorem \ref{tauthm} gives an oracle inequality for the PLSE built on $\tau$-mixing sequences. This inequality is not pathwise and the constants involved in the penalty term are not optimal. This is due to technical reasons, mainly because we use the coupling result (\ref{coup2}) instead of (\ref{coup1}). However, we recover the same kind of oracle inequality as in the i.i.d. framework (Birg\'e and Massart \cite{BM97}) under weak assumptions on the mixing coefficients since we only require arithmetical {\bf {[AR]}} $\tau$-mixing assumptions on the process $(X_n)_{n\in\Z}$. This is the first result for these processes up to our knowledge.\\
Let us mention here Theorem 4.1 in Comte $\&$ Merlev\`ede \cite{CM02}. They consider $\alpha$-mixing processes (for a definition of the coefficient $\alpha$ and its properties, we refer to Rio \cite{Ri00}). They make geometrical {\bf {[GEO]}} $\alpha$-mixing assumptions on the processes and consider penalties of order $L\log(n)D_m/n$ to get an oracle inequality. This leads to a logarithmic loss in the rates of convergence. They get the optimal rate under an extra assumption (namely Assumption $[\Lip]$ in Section 3.2). There exist random processes that are $\tau$-mixing and not $\alpha$-mixing (see Dedecker $\&$ Prieur \cite{DP05}), however, the comparison of these coefficients is difficult in general and our method can not be applied in this context. \\
The constants $c_0,\gamma_1,n_o$ are given in the end of the proof.\\
{\bf Remark :}
Inequality (2.6) can be improved under stronger assumptions on $s$. For example, when $s$ is bounded, we have $\tilde{\beta}_k\leq C\sqrt{\tau_{k}}$. Under this assumption and $\theta>3$, we can prove that the estimator $\tilde{s}$ satisfies the inequality 
\begin{equation*}
\E\left(\|\tilde{s}-s\|_2^2\right)\leq (1+\epsilon_n)\left(\inf_{m\in \M_n,\;D_m \geq c_0(\log n)^{\gamma_1}}\|s-s_m\|_2^2+\pen(m)\right)+\frac{(\log n)^{\kappa(\theta+1)}}{n^{(\theta-3)/2}}.
\end{equation*}
When $\theta<5$, the extra term $(\log n)^{\kappa(\theta+1)}/n^{(\theta-3)/2}$ may be larger than the main term $\inf_{m\in \M_n,\;D_m \geq c_0(\log n)^{\gamma_1}}\|s-s_m\|_2^2+\pen(m)$. In this case, we don't know if our control remains optimal. On the other hand, Proposition \ref{minimax} ensures that $\tilde{s}$ is adaptive over the class of Besov balls when $\theta\geq 5$.
\section{Minimax results}
\subsection{Approximation results on Besov spaces}
{\bf{Besov balls.}}\\
Throughout this section, $\Lambda=\{(j,k),\;j\in \N,\; k\in \Z\}$ and $\{\psi_{j,k},\;(j,k)\in \Lambda\}$ denotes an $r$-regular wavelet basis as introduced in Section 4.1. Let $\alpha,p$ be two positive numbers such that $\alpha+1/2-1/p>0$. For all functions $t\in L_2(\mu)$, $t=\sum_{(j,k)\in \Lambda} t_{j,k}\psi_{j,k}$, we say that $t$ belongs to the Besov ball $B_{\alpha,p,\infty}(M_1)$ on the real line if $\left\|t\right\|_{\alpha,p,\infty}\leq M_1$ where 
$$\left\|t\right\|_{\alpha,p,\infty}=\sup_{j\in\N}2^{j(\alpha+1/2-1/p)}\left(\sum_{k\in\Z}|t_{j,k}|^p\right)^{1/p}.$$ 
It is easy to check that if $p\geq 2$ $B_{\alpha,p,\infty}(M_1) \subset B_{\alpha,2,\infty}(M_1) $ so that upper bounds on $B_{\alpha,2,\infty}(M_1)$ yield upper bounds on $B_{\alpha,p,\infty}(M_1)$.\\
{\bf{Approximation results on Besov spaces.}}\\
We have the following result (Birg\'e $\&$ Massart \cite{BM97} Section 4.7.1). Suppose that the support of $s$ equals $[0,1]$ and that $s$ belongs to the Besov ball $B_{\alpha,2,\infty}(1)$, then whenever $r>\alpha-1$,
\begin{equation}\label{ApproxN}
\left\|s-s_m\right\|_2^2\leq \frac{\left\|s\right\|_{\alpha,2,\infty}^2}{4(4^{\alpha}-1)}2^{-2J_m\alpha}\leq\frac{(2A)^{2\alpha}\left\|s\right\|_{\alpha,2,\infty}^2}{4(4^{\alpha}-1)}D_m^{-2\alpha}
\end{equation}
\subsection{Minimax rates of convergence for the PLSE}
We can derive from Theorems \ref{betathm} and \ref{tauthm} adaptation results to unknown smoothness over Besov Balls.
\begin{prop}\label{minimaxb}
Assume that the process $(X_n)_{n\in\Z}$ is stricly stationary and arithmetically {\bf{[AR]}} $\beta$-mixing with mixing rate $\theta>2$ and that its marginal distribution admits a density $s$ with respect to the Lebesgue measure $\mu$, that $s$ is supported in $[0,1]$ and that $s\in L^2(\mu)$. For all $\alpha,M_1>0$, the PLSE $\tilde{s}$ defined in Theorem \ref{betathm} for the collection of models {\bf{[W]}} satisfies
$$\forall \kappa>2,\;\sup_{s\in B_{\alpha,2,\infty}(M_1)}\p\left(\left\|\tilde{s}-s\right\|_2^2> L_{M_1,\alpha,\theta}n^{-2\alpha/(2\alpha+1)}\right)\leq \frac{L_{M_1}(\log n)^{(\theta+2)\kappa}}{n^{\theta/2}}.$$
\end{prop}

\begin{prop}\label{minimax}
Assume that the process $(X_n)_{n\in\Z}$ is stricly stationary and arithmetically {\bf{[AR]}} $\tau$-mixing with mixing rate $\theta>5$ and that its marginal distribution admits a density $s$ with respect to the Lebesgue measure $\mu$, that $s$ is supported in $[0,1]$ and that $s\in L^2(\mu)$. For all $\alpha,M_1>0$, the PLSE $\tilde{s}$ defined in Theorem \ref{tauthm} satisfies
$$\sup_{s\in B_{\alpha,2,\infty}(M_1)}\E\left(\left\|\tilde{s}-s\right\|_2^2\right)\leq L_{M_1,\alpha,\theta}n^{-2\alpha/(2\alpha+1)}.$$
\end{prop}
{\bf Remark:}
Proposition \ref{minimax} can be compared to Theorem 3.1 in Gannaz $\&$ Wintenberger \cite{GW09}. They prove near minimax results for the thresholded wavelet estimator introduced by Donoho {\it{et al.}} \cite{DJKP} in a $\tilde{\phi}$-dependent setting (for a definition of the coefficient $\tilde{\phi}$, we refer to Dedecker $\&$ Prieur \cite{DP05}). Basically, with our notations, their result can be stated as follows: if $(X_n)_{n\in\Z}$ is $\tilde{\phi}$-mixing with $\tilde{\phi}_1(r)\leq Ce^{-ar^b}$ for some constants $C,a,b$, then the thresholded wavelet estimator $\hat{s}$ of $s$ satisfies 
$$\forall \alpha>0,\;\forall p>1,\;\sup_{s\in B_{\alpha,p,\infty}(M_1)\cap L^{\infty}(M)}\E\left(\left\|\hat{s}-s\right\|_2^2\right)\leq L_{M,M_1,\alpha,p}\left(\frac{\log n}{n}\right)^{2\alpha/(2\alpha+1)}.$$
The main advantage of their result is that they can deal with Besov balls with regularity $1<p<2$. However, in the regular case, when $p\geq 2$, we have been able to remove the extra $\log n$ factor. Moreover, our result only requires arithmetical {\bf{[AR]}} rates of convergence for the mixing coefficients and we do not have to suppose that $s$ is bounded.

\section{Proofs.}
\subsection{Proofs of the minimax results.}
\begin{proof}{of Proposition \ref{minimaxb}:}
Let $\alpha>0$ and $M_1>0$ and assume that $s\in B_{\alpha,2,\infty}(M_1)$. Let $\tilde{\M}_n=\{m\in\M_n,D_m>c_0(\log n)^{\gamma_1}\}$. By Theorem \ref{betathm}, there exists a constant $L_{\theta}>0$ such that
\begin{equation}\label{min1}
\p\left(\left\|\tilde{s}-s\right\|_2^2>L_{\theta}\inf_{m\in\tilde{\M}_n}\left\{ \left\|s-s_{m}\right\|_2^2+\frac{D_m}{n}\right\}\right)\leq \frac{L_s(\log n)^{(\theta+2)\kappa}}{n^{\theta/2}}.
\end{equation}
It appears from the proof of Theorem \ref{betathm} that the constant $L_s$ depends only on $\|s\|_2$ and that it is a nondecreasing function of $\|s\|_2$ so that $L_s$ can be uniformly bounded over $B_{\alpha,2,\infty}(M_1)$ by a constant $L_{M_1}$ so that, by (\ref{min1})  
\begin{equation*}
\p\left(\left\|\tilde{s}-s\right\|_2^2>L_{\theta}\inf_{m\in\tilde{\M}_n}\left\{ \left\|s-s_m\right\|_2^2+\frac{D_m}{n}\right\}\right)\leq \frac{L_{M_1}(\log n)^{(\theta+2)\kappa}}{n^{\theta/2}}.
\end{equation*}
In particular, for a model $m$ in $\M_n$ with dimension $D_m$ such that 
$$c_0(\log n)^{\gamma_1}\leq L_1n^{1/(2\alpha+1)}\leq D_m \leq L_2 n^{1/(2\alpha+1)},$$
we have
\begin{equation*}
\p\left(\left\|\tilde{s}-s\right\|_2^2>L_{\theta} \left(\left\|s-s_m\right\|_2^2+\frac{D_m}{n}\right)\right)\leq \frac{L_{M_1}(\log n)^{(\theta+2)\kappa}}{n^{\theta/2}}.
\end{equation*}
Since $s$ belongs to $B_{\alpha,2,\infty}(M_1)$, we can use Inequality (\ref{ApproxN}) to get 
$$\left\|s-s_m\right\|_2^2\leq L_{\alpha,M_1}D_m^{-2\alpha}.$$
Thus we obtain
$$\p\left(\left\|\tilde{s}-s\right\|_2^2> L_{M_1,\alpha,\theta}n^{-2\alpha/(2\alpha+1)}\right)\leq \frac{L_{M_1}(\log n)^{(\theta+2)\kappa}}{n^{\theta/2}}.\square$$
\end{proof}
\begin{proof}{of Proposition \ref{minimax}:}
Let $\alpha>0$ and $M_1>0$ and assume that $s\in B_{\alpha,2,\infty}(M_1)$. By Theorem \ref{tauthm}, we have
\begin{eqnarray*}
\E\left(\left\|\tilde{s}-s\right\|_2^2\right)&\leq& L_{\theta} \left(\inf_{m\in\tilde{\M}_n}\lbrace\left\|s-s_{m}\right\|_2^2+\frac{D_m}{n}\rbrace\right).\\
\end{eqnarray*}
Inequality (\ref{ApproxN}) leads to $\left\|s-s_{m}\right\|_2^2\leq L_{\alpha,M_1}D_m^{-2\alpha},$ so that for a model $m$ in $\tilde{M}_n$ with dimension $D_m$ such that 
$$c_0(\log n)^{\gamma_1}\leq L_1n^{1/(2\alpha+1)}\leq D_m \leq L_2 n^{1/(2\alpha+1)},$$ 
we find 
$$\E\left(\left\|\tilde{s}-s\right\|_2^2\right)\leq L_{\theta,\alpha,M_1}n^{-2\alpha/(2\alpha+1)}.\square$$
\end{proof}
\subsection{Proof of Theorem \ref{betathm}:}
For all $m_o$ in $\M_n$, we have, by definition of $\hat{m}$
\begin{eqnarray*}
\gamma_n(\tilde{s})+\pen(\hat{m})&\leq &\gamma_n(\hat{s}_{m_o})+\pen(m_o)\\
P\gamma(\tilde{s})+\nu_n\gamma(\tilde{s})+\pen(\hat{m})&\leq&P\gamma(\hat{s}_{m_o})+\nu_n\gamma(\hat{s}_{m_o})+\pen(m_o)\\
P\gamma(\tilde{s})-P\gamma(s)-2\nu_n\tilde{s}+\pen(\hat{m})&\leq&P\gamma(\hat{s}_{m_o})-P\gamma(s)-2\nu_n\hat{s}_{m_o}+\pen(m_o)
\end{eqnarray*}
Since for all $t\in L_2(\mu)$, $P\gamma(t)-P\gamma(s)=\|t-s\|_2^2$, we have 
\begin{equation}\label{dec}
\left\|s-\tilde{s}\right\|_2^2\leq\left\|s-\hat{s}_{m_o}\right\|_2^2+\pen(m_o)-V(m_o)-(\pen(\hat{m})-V(\hat{m}))-2\nu_n(s_{m_o}-s_{\hat{m}}),
\end{equation}
where, for all $m\in\M_n$
$$V(m)=2\nu_n(\hat{s}_m-s_m)=2\sum_{(j,k)\in m}\nu_n^2(\psi_{j,k}).$$
This decomposition is different from the one used in Birg\'e $\&$ Massart \cite{BM97} and in Comte $\&$ Merlev\`ede \cite{CM02}. It allows to improve the constant in the oracle inequality in the $\beta$-mixing case. Moreover, we choose to prove an oracle inequality of the form (\ref{oracleP}) for $\beta$-mixing sequences, which allows to assume only $\theta>2$ instead of $\theta>3$. 
Let us now give a sketch of the proof:
\begin{enumerate}
\item we build an event $\Omega_C$ with $\p(\Omega^c_C)\leq p\beta_q$ such that, on $\Omega_C$, $\nu_n=\nu_n^{*}$, where $\nu_n^{*}$ is built with independent data. A suitable choice of the integers $p$ and $q$ leads to $p\beta_q\leq C(\ln n)^{r}n^{-\theta/2}$.
\item We use the concentration's inequality (\ref{Massart}) of Birg\'e $\&$ Massart \cite{BM97} for $\chi^2$-type statistics, derived from Talagrand's inequality. This allows us to find $p_1(m)$ such that on an event $\Omega_1$ with $\p(\Omega_1^c\cap \Omega_C)\leq L_{1,s}c_n$
$$\sup_{m\in\M_n}\left\{V(m)-p_1(m)\right\}\leq 0.$$
$c_n<C(\ln n)^rn^{-\theta/2}$ and $L_{1,s}$ is some constant depending on $s$.
\item From Bernstein's inequality, we prove that, for all $m,m'\in \M_n$, there exists $p_2(m,m')$ such that, for all $\eta>0$, on an event $\Omega_2$ with $\p(\Omega^c_2\cap\Omega_C)\leq L_{2,s}c_n,$
\begin{eqnarray*}
\sup_{m,m'\in\M_n}\left\{\nu_n(s_m-s_{m'})-\frac{\eta}{2} p_2(m,m')-\frac{\left\|s_m-s_{m'}\right\|_2^2}{2\eta}\right\}\leq 0.
\end{eqnarray*}
Moreover, for all $m,m'\in\M_n$, $p_2(m,m')\leq p_2(m,m)+p_2(m',m').$
\item We have $\left\|s_{\hat{m}}-s_{m_o}\right\|_2^2\leq  \left\|s_{\hat{m}}-s\right\|_2^2+\left\|s-s_{m_o}\right\|_2^2$ because $s_{\hat{m}}-s_{m_o}$ is either the projection of $s_{\hat{m}}-s$ onto $S_{m_o}$ or the projection of $s-s_{m_o}$ onto $S_{\hat{m}}$. Take $\pen(m)\geq p_1(m)+\eta p_2(m,m)$, we have, on $\Omega_1\cap\Omega_2\cap\Omega_C$
\end{enumerate}
\begin{eqnarray}
\left\|s-\tilde{s}\right\|_2^2&\leq&\left\|s-\hat{s}_{m_o}\right\|_2^2-\frac{V_{m_o}}2+\pen(m_o)-\frac{V_{m_o}}2\\
&&-(\pen(\hat{m})-p_1(\hat{m}))-(p_1(\hat{m})-V(\hat{m}))-2\nu_n(s_{m_o}-s_{\hat{m}})\nonumber\\
&\leq&\left\|s-s_{m_o}\right\|_2^2+\pen(m_o)-\frac{V(m_o)}2-\eta p_2(\hat{m},\hat{m})\nonumber\\
&&+\eta p_2(\hat{m},m_o)+\frac{\|s_{m_o}-s_{\hat{m}}\|_2^2}{\eta}\label{afin}\\
\left(1-\frac{1}{\eta}\right)\left\|s-\tilde{s}\right\|_2^2&\leq&(1+\frac{1}{\eta})\left\|s-s_{m_o}\right\|_2^2+\pen(m_o)+\eta p_2(m_o,m_o).\label{fin}
\end{eqnarray}
In (\ref{afin}), we used that $V(m_o)=2\|s_{m_o}-\hat{s}_{m_o}\|_2^2\geq 0$. In (\ref{fin}), we used that $V_{m_o}\geq 0$. Pythagoras Theorem gives 
$$\left\|s-\hat{s}_{m_o}\right\|_2^2-\frac{V(m_o)}2=\left\|s-s_{m_o}\right\|_2^2\;\textrm{and}_; \|s-s_{\hat{m}}\|_2^2\leq \left\|s-\tilde{s}\right\|_2^2.$$
Finally, we prove that we can choose $\eta=(\log n)^{\gamma}$, with $\gamma>0$ such that $\eta p_2(m_o,m_o)=o(\pen(m_o))$ and we conclude the proof of (\ref{betathm}) from the previous inequalities.\\
We decompose the proof in several claims corresponding to the previous steps.\\
{\bf{Claim 1 :}} For all $l=0,...,p-1$, let us define $A_l=(X_{2lq+1},\; ...,X_{(2l+1)q})$ and $B_l=(X_{(2l+1)q+1},...,X_{(2l+2)q})$. There exist random vectors $A_l^*=(X_{2lq+1}^*,...,X_{(2l+1)q}^*)$ and $B_l^*=(X_{(2l+1)q+1}^*,...,X_{(2l+2)q}^*)$ such that for all $l=0,...,p-1$ :
\begin{enumerate}
\item $A_l^*$ and $A_l$ have the same law, 
\item $A_l^*$ is independent of $A_0,...,A_{l-1},A_0^*...,A_{l-1}^*$ 
\item $\p(A_l\neq A_l^*)\leq \beta_q$ 
\end{enumerate}
the same being true for the variables $B_l$.
\begin{proof}{of Claim 1 :}
The proof is derived from Berbee's lemma, we refer to Proposition 5.1 in Viennet \cite{Vi97} for further details about this construction.$\square$
\end{proof}
Hereafter, we assume that, for some $\kappa>2$, $\sqrt{n}(\log n)^{\kappa}/2\leq p\leq \sqrt{n}(\log n)^{\kappa}$ and for the sake of simplicity that $pq=n/2$, the modifications needed to handle the extra term when $q=[n/(2p)]$ being straightforward. Let $\Omega_C=\left\lbrace  \forall l=0,...,p-1\; A_l=A_l^*,\;B_l=B_l^*\right\rbrace$. We have $$\p(\Omega_C^c)\leq 2p\beta_q\leq 2^{2+\theta}\frac{(\log n)^{(\theta+2)\kappa}}{n^{\theta/2}}.$$ 
Let us first deal with the quadratic term $V(m)$.\\
{\bf{Claim 2 :}} {\it{Under the assumptions of Theorem \ref{betathm}, let $\epsilon>0$, $1<\gamma<\kappa/2$. We define $L_1^2= 2\Phi^2\kappa_1$, $L_2^2=8\Phi^{3/2}\sqrt{\kappa_2}$, $L_3=2\Phi\kappa(\epsilon)$ and
\begin{equation}\label{L1m}
L_{1,m}=4\left((1+\epsilon)L_1+L_2\sqrt{\frac{(\log n)^{\gamma}}{D_m^{1/4}}}+\frac{L_3}{(\log n)^{\kappa-\gamma}}\right)^2.
\end{equation}
Then, we have 
\begin{equation*}
\p\left(\sup_{m\in\M_n}\left\{V(m)- \frac{L_{1,m}D_m}{n}\right\}\geq 0\cap \Omega_C\right)\leq L_{s,\gamma}\exp\left(-\frac{(\log n)^{\gamma}}{\sqrt{\|s\|_2}}\right).
\end{equation*}
where $L_{s,\gamma}=2\sum_{D=1}^{\infty} \exp(-(\log D)^{\gamma}/\left\|s\right\|_2^{1/2}).$ In particular, for all $r>0$, there exists a constant $L'_{s,r}$ depending on $\left\|s\right\|_2$, such that 
$$\p\left(\sup_{m\in\M_n}\left\{V(m)- \frac{L_{1,m}D_m}{n}\right\}\geq 0\cap \Omega_C\right)\leq\frac{L'_{s,r}}{n^r}.$$}}
{\bf Remark :}\label{C2}
When $(L_2/L_1)^8(\log n)^{4(2\kappa-\gamma)}\leq D_m\leq n$, we have 
$$L_{1,m}\leq \left[1+\epsilon+\left(1+\frac{\sqrt{2}\kappa(\epsilon)}{\sqrt{\kappa_1}}\right)(\log n)^{-(\kappa-\gamma)}\right]^24L_1^2.$$
\begin{proof}{of Claim 2 :}
Let $P_n^*(t)=\sum_{i=1}^nt(X_i^*)/n$ and $\nu_n^*(t)=(P_n^*-P)t$, we have
$$V(m){\bf{1}}_{\Omega_C}=2\sum_{(j,k)\in m}(\nu_n^*)^2(\psi_{j,k}){\bf{1}}_{\Omega_C}.$$
Let $B_1(S_m)=\left\{t\in S_m;\;\left\|t\right\|_2\leq 1\right\}$. $\forall t\in B_1(S_m)$, let $\bar{t}(x_1,...,x_q)=\sum_{i=1}^q t(x_i)/2q$ and for all functions $g:\R^q\rightarrow \R$ let 
$$P^*_{A,p}g=\frac{1}{p}\sum_{j=0}^{p-1}g(A^*_j),\;P^*_{B,p}g=\frac{1}{p}\sum_{j=0}^{p-1}g(B^*_j),\;\bar{P}g=\int g\p_A(d\mu),$$
$$\textrm{and}\;\bar{\nu}_{A,p}g=(P^*_{A,p}-\bar{P})g,\;\bar{\nu}_{B,p}g=(P^*_{B,p}-\bar{P})g.$$
Now we have
$$\sum_{(j,k)\in m}(\nu_n^*)^2(\psi_{j,k})\leq 2\sum_{(j,k)\in m}\bar{\nu}_{A,p}^2\bar{\psi}_{j,k}+2\sum_{(j,k)\in m}\bar{\nu}_{B,p}^2\bar{\psi}_{j,k}.$$
In order to handle these terms, we use Proposition \ref{Massart} which is stated in Section 7. Taking
$$B_m^2=\sum_{(j,k)\in m}\Var(\bar{\psi}_{j,k}(A_1)),\;V_m^2=\sup_{t\in B_1(S_m)}\Var(\bar{t}(A_1)),\;\textrm{and}\;H_m^2=\left\|\sum_{(j,k)\in m}(\bar{\psi}_{j,k})^2\right\|_{\infty},$$
we have
\begin{equation}\label{concb}
\forall x>0,\;\p\left(\sqrt{\sum_{(j,k)\in m}\bar{\nu}_{A,p}^2\bar{\psi}_{j,k}}\geq \frac{(1+\epsilon)}{\sqrt{p}}B_m+V_m\sqrt{\frac{2x}{p}}+\kappa(\epsilon)\frac{H_m x}{p}\right)\leq e^{-x}.
\end{equation}
In order to evaluate $B_m$, $V_m$ and $H_m$, we use Viennet's inequality (\ref{Viennet}). There exists a function $b$ such that, for all $p= 1,2$, $P|b|^p\leq \kappa_p$ where $\kappa_p$ is defined in (\ref{1}) and for all functions $t\in L_2(\bar{P})$, 
$$\Var(\bar{t}(A_1))\leq \frac{1}{q}Pbt^2.$$
Thus 
\begin{equation*}
B_m^2=\sum_{(j,k)\in m}\Var(\bar{\psi}_{j,k}(A_1))\leq \frac{1}{q}\sum_{(j,k)\in m}Pb\psi^2_{j,k}\leq \left\|\sum_{(j,k)\in m}\psi^2_{j,k}\right\|_{\infty}\frac{\kappa_1}{q}.
\end{equation*}
From Assumption {\bf {[$M_2$]}}, $\left\|\sum_{(j,k)\in m}\psi^2_{j,k}\right\|_{\infty}\leq \Phi^2D_m,$ thus,
\begin{equation}\label{bm}
B_m^2\leq \frac{\Phi^2\kappa_1 D_m}{q}.
\end{equation}
From Viennet's and Cauchy-Schwarz inequalities 
$$V_m^2=\sup_{t\in B_1(S_m)}\Var(\bar{t}(A_1))\leq \sup_{t\in B_1(S_m)}\frac{Pbt^2}{q}\leq \sup_{t\in B_1(S_m)}\left\|t\right\|_{\infty}\frac{(Pt^2)^{1/2}(Pb^2)^{1/2}}{q}.$$
Since $t\in B_1(S_m)$, we have by Cauchy-Schwarz inequality 
$$(Pt^2)^{1/2}\leq (\left\|t\right\|_{\infty}\left\|t\right\|_2\left\|s\right\|_2)^{1/2}\leq (\left\|t\right\|_{\infty}\left\|s\right\|_2)^{1/2}.$$
From Assumption {\bf {[$M_2$]}}, we have $\left\|t\right\|_{\infty}\leq \Phi\sqrt{D_m}$, and from Viennet's inequality $Pb^2\leq\kappa_2<\infty$, thus we obtain
\begin{equation}\label{vm}
V_m^2\leq \Phi^{3/2}(\left\|s\right\|_2\kappa_2)^{1/2}\frac{D_m^{3/4}}{q}.
\end{equation}
Finally, from Assumption {\bf {[$M_2$]}}, we have, using Cauchy-Schwarz inequality
\begin{equation}\label{hm}
H_m^2=\left\|\sum_{(j,k)\in m}\bar{\psi}_{j,k}^2\right\|_{\infty}\leq \frac{1}{4}\left\|\sum_{(j,k)\in m}\psi_{j,k}^2\right\|_{\infty}\leq \frac{\Phi^2D_m}{4}.
\end{equation}
Let $y_n>0$. We define 
$$L_m=\left((1+\epsilon)L_1+L_2\sqrt{\frac{(\log D_m)^{\gamma}+y_n}{2D_m^{1/4}}}+L_3\frac{(\log D_m)^{\gamma}+y_n}{2(\log n)^{\kappa}}\right)^2.$$
We apply Inequality (\ref{concb}) with $x=((\log D_m)^{\gamma}+ y_n)/\left\|s\right\|_2^{1/2}$ and the evaluations (\ref{bm}), (\ref{vm}) and (\ref{hm}). Recalling that $1/p\leq 2/(\sqrt{n}(\log n)^{\kappa})$, this leads to 
\begin{equation*}
\p\left(\sum_{(j,k)\in m}\bar{\nu}_{A,p}^2\bar{\psi}_{j,k}\geq \frac{L_mD_m}{n}\right)\leq \exp\left(-\frac{(\log D_m)^{\gamma}}{\sqrt{\left\|s\right\|_2}}\right)\exp(-\frac{y_n}{\sqrt{\|s\|_2}}).
\end{equation*}
In order to give an upper bound on $H_mx$, we used that the support of $s$ in included in $[0,1]$, thus 
$$1=\left\|s\right\|_1\leq \left\|s\right\|_2.$$
The result follows by taking $y_n=(\log n)^{\gamma}\geq (\log D_m)^{\gamma}.\square$
\end{proof}
{\bf{Claim 3.}} {\it{We keep the notations $\kappa/2>\gamma>1$, $L_2$ of the proof of Claim 2. For all $m,m'\in \M_n$ we take 
\begin{equation}\label{Lmm}
L_{m,m'}=4\left(L_2\sqrt{\frac{(\log n)^{\gamma}}{(D_m\vee D_{m'})^{1/4}}}+\frac{4\Phi}{3(\log n)^{\kappa-\gamma}}\right)^2,
\end{equation}
we have, for all $\eta>0$,
\begin{equation*}
\p\left(\sup_{m,m'\in\M_n}\nu^*_{n}(s_m-s_{m'})-\frac{\left\|s_m-s_{m'}\right\|_2^2}{2\eta}-\frac{\eta}{2}\frac{L_{m,m'}(D_m\vee D_{m'})}{n}>0\right)\leq L_{s,\gamma}e^{-\frac{(\log n)^{\gamma}}{\left\|s\right\|_2^{1/2}}}
\end{equation*}
with $L_{s,\gamma}=2 \sum_{m,m'\in\M_n} e^{-\frac{(\log (D_m\vee D_{m'}))^{\gamma}}{\left\|s\right\|_2^{1/2}}}.$}}\\
{\bf Remark :}\label{RC3}
The constant $L_{s,\gamma}$ is finite since for all $x,y >0$, $(\log (x\vee y))^{\gamma}\geq ((\log x)^{\gamma}+(\log y)^{\gamma})/2$.\\
As in Claim 2, when $(L_2/L_1)^8(\log n)^{4(2\kappa-\gamma)}\leq D_m\leq n$, we have 
$$L_{m,m'}\leq \left(1+\frac{2^{3/2}}{3\sqrt{\kappa_1}}\right)^2(\log n)^{-2(\kappa-2\gamma)}4L_1^2.$$

\begin{proof}{of Claim 3.}
We keep the notations of the proof of Claim 2 and for $m,m'\in \M_n$, let $t_{m,m'}=(s_m-s_{m'})/\left\|s_m-s_{m'}\right\|_2$. We use the inequality $2ab\leq a^2\eta^{-1}+b^2\eta$, which holds for all $a,b\in\R$, $\eta>0$. This leads to 
\begin{eqnarray*}
\nu^*_n(s_m-s_{m'})&=&\left\|s_m-s_{m'}\right\|_2\nu^*_n(t_{m,m'})\leq\frac{\left\|s_m-s_{m'}\right\|_2^2}{2\eta}+\frac{\eta}{2}\left(\nu^*_n(t_{m,m'})\right)^2\\
&=&\frac{\left\|s_m-s_{m'}\right\|_2^2}{2\eta}+\frac{\eta}{2}\left(\bar{\nu}_{A,p}(\bar{t}_{m,m'})+\bar{\nu}_{B,p}(\bar{t}_{m,m'})\right)^2\\
&\leq & \frac{\left\|s_m-s_{m'}\right\|_2^2}{2\eta}+ \eta ( \bar{\nu}_{A,p}(\bar{t}_{m,m'}))^2+ \eta ( \bar{\nu}_{B,p}(\bar{t}_{m,m'}))^2.
\end{eqnarray*}
Now from Bernstein's inequality (see Section 7), we have
\begin{equation}\label{ber}
\forall x>0,\;\p\left(\bar{\nu}_{A,p}(\bar{t}_{m,m'})>\sqrt{\frac{2\Var(\bar{t}_{m,m'}(A_1))x}{p}}+\frac{\Vert\bar{t}_{m,m'}\Vert_{\infty}x}{3p}\right)\leq e^{-x}.
\end{equation}
From Viennet's and Cauchy-Schwarz inequalities, we have 
$$\Var(\bar{t}_{m,m'}(A_1))\leq \frac{Pbt_{m,m'}^2}{q}\leq \frac{\Vert t_{m,m'}\Vert_{\infty}\sqrt{Pb^2Pt_{m,m'}^2}}{q}.$$
Moreover
$$Pb^2\leq \kappa_2,\;Pt_{m,m'}^2\leq \Vert t_{m,m'}\Vert_{\infty}\Vert t_{m,m'}\Vert_2\Vert s\Vert_2.$$
Since $t_{m,m'}\in S_m\cup S_{m'}$ and $\Vert t_{m,m'}\Vert_2 = 1$, we have, from Assumption {\bf {[$M_2$]}} $\Vert t_{m,m'}\Vert_{\infty}\leq \Phi\sqrt{D_m\vee D_{m'}}$. Let $y_n>0$. We apply Inequality (\ref{ber}) with  $x=[(\log (D_m\vee D_{m'}))^{\gamma}+y_n]/\left\|s\right\|_2^{1/2}$. We define 
$$\frac{L'_{m,m'}}4 = \left(L_2\sqrt{\frac{(\log (D_m\vee D_{m'}))^{\gamma}+y_n}{2(D_m\vee D_{m'})^{1/4}}}+\frac{4\Phi\left[(\log (D_m\vee D_{m'}))^{\gamma}+y_n\right]}{6(\log n)^{\kappa}}\right)^2,$$
we have
\begin{equation*}
\p\left(\bar{\nu}_{A,p}(\bar{t}_{m,m'})>\sqrt{\frac{L'_{m,m'}(D_m\vee D_{m'})}{4n}}\right)\leq \exp\left(-\frac{(\log(D_m\vee D_{m'}))^{\gamma}}{\left\|s\right\|_2^{1/2}}\right)e^{-y_n/\left\|s\right\|_2^{1/2}}.
\end{equation*}
The result follows by taking $y_n=(\log n)^{\gamma}$ and using $2\leq D_m\leq n.$\\
{\bf{Conclusion of the proof:}}\\
Let $\eta>0$ and $\pen'(m)\geq(L_{1,m}+\eta L_{m,m})D_m/n$ where $L_{1,m}$ and $L_{m,m}$ are defined respectively by (\ref{L1m}) and (\ref{Lmm}). From Claims 1, 2 and 3 and (\ref{fin}), we obtain that, for all $m_o$ and with probability larger than $L_{s,\theta}(\log n)^{(\theta+2)\kappa}n^{-\theta/2}$
\begin{equation}\label{fp}
(1-\frac{1}{\eta})\left\|s-\tilde{s}\right\|_2^2\leq (1+\frac{1}{\eta})\left\|s-s_{m_o}\right\|_2^2+\pen'(m_o)+\eta L(m_o,m_o)\frac{D_{m_o}}{n}.
\end{equation}
Assume that $D_{m}\geq (L_2/L_1)^8(\log n)^{4(2\kappa-\gamma)}$, then we have from remarks \ref{C2} and \ref{RC3}
\begin{eqnarray*}
L_{1,m}&\leq& \left[1+\epsilon+\left(1+\frac{2\kappa(\epsilon)}{\sqrt{\kappa_1}}\right)(\log n)^{-(\kappa-2\gamma)}\right]^24L^2_1\;\textrm{and}\\
L_{m,m}&\leq&\left(1+\frac{2^{3/2}}{3\sqrt{\kappa_1}}\right)^2(\log n)^{-2(\kappa-\gamma)}4L^2_1.
\end{eqnarray*}
Take $\eta=(\log n)^{\kappa-\gamma}$, we have $(L_{1,m_o}+\eta L_{m_o,m_o})D_{m_o}/n\leq C\pen(m_o)$. Fix $\epsilon>0$ such that $[1+\epsilon]^2<K/4$. Since $\kappa>\gamma$, for $n\geq n_o$, we have $L_{1,m}+\eta L_{m,m}\leq KL_1^2$, thus, inequality (\ref{Beta}) follows follows from (\ref{fp}) as soon as $n>n_o$. We remove the condition $n>n_o$ by improving the constant $L_s$ in (\ref{Beta}) if necessary.$\square$
\end{proof}

\subsection{Proof of Theorem \ref{tauthm}.}
The proof follows the previous one, the main difference is that the coupling lemma (Claim 1) as well as the covariance inequalities are much harder to handle in the $\tau$-mixing case. This leads to more technical computations to recover the results obtained in the $\beta$-mixing case (see Claims 2, 3 and the proof of inequality (\ref{rb})). We start with the decomposition (\ref{dec}). As in the previous proof, the decomposition of the risk given in Birg\'e $\&$ Massart \cite{BM97} or in Comte $\&$ Merlev\`ede \cite{CM02} could be used. This leads to a loss in the constant in front of the main term in (\ref{tau}) without avoiding any of the main difficulties. We divide the proof in four claims. \\ 
{\bf{Claim 1 :}} {\it{For all $l=0,...,p-1$, let us denote by $A_l=(X_{2lq+1},...,X_{(2l+1)q})$ and $B_l=(X_{(2l+1)q+1},...,X_{(2l+2)q})$. There exist random vectors $A_l^*=(X_{2lq+1}^*,...,X_{(2l+1)q}^*)$ and $B_l^*=(X_{(2l+1)q+1}^*,...,X_{(2l+2)q}^*)$ such that for all $l=0,...,p-1$ :
\begin{itemize}
\item $A_l^*$ and $A_l$ have the same law, 
\item $A_l^*$ is independent of $A_0,...,A_{l-1},A_0^*...,A_{l-1}^*$ 
\item $\E(|A_l-A_l^*|_{q})\leq q\tau_q$ 
\end{itemize}
the same being true for the variables $B_l$.}}
\begin{proof}{of Claim 1 :} We use the same recursive construction as Viennet \cite{Vi97}. \\
Let $(\delta_j)_{0\leq j\leq p-1}$ be a sequence of independent random variables uniformly distributed over $[0,1]$ and independent of the sequence $(A_j)_{0\leq j\leq p-1}$.
Let $A_0^*=(X_{1}^*,...,X_q^*)$ be the random variable given by equality (\ref{coup2}) for $\M=\sigma(X_i,\;i\leq-q)$, $A_0$ and $\delta_0$.\\
Now suppose that we have built the variables $A_l^*$ for $l< l'$. From equality (\ref{coup2}) applied to the $\sigma$-algebra $\sigma(A_l,A_l^*,\; l<l')$, $A_{l'}$ and $\delta_{l'}$, there exists a random variable $A_{l'}^*$ satisfying the hypotheses of Claim 1. \\
We build in the same way the variables $B_l^*$ for all $l=0,...,p-1$. $\square$
\end{proof}
We keep the notations $\nu_n^*,\bar{\nu}_{A,p},\bar{\nu}_{B,p}$, $\bar{t}$ and $B_1(S_m)$ that we introduced in the proof of Theorem \ref{betathm}. As in the proof of Theorem \ref{betathm}, we assume that, for some $\kappa>2$, $\sqrt{n}(\log n)^{\kappa}/2\leq p\leq \sqrt{n}(\log n)^{\kappa}$ and for the sake of simplicity that $pq=n/2$, the modifications needed to handle the extra term when $q=[n/(2p)]$ being straightforward. We have
\begin{eqnarray}
V(\hat{m})&=&\sum_{(j,k)\in\hat{m}}\nu_n^2(\psi_{j,k})\leq 2\sum_{(j,k)\in\hat{m}}(P_n-P_n^*)^2(\psi_{j,k})+2\sum_{(j,k)\in\hat{m}}(\nu_n^*)^2(\psi_{j,k})\label{EP}
\end{eqnarray}
{\bf{Claim 2 :}} {\it{There exists a constant $L=L_{A,K_L,K_{\infty},\kappa,\theta}$ such that
\begin{equation}\label{reste}
\E\left(\sum_{j,k\in\hat{m}}\left((P_n-P_n^*)(\psi_{j,k})\right)^2\right)\leq L\frac{(\log n)^{\kappa(\theta+1)}}{n^{(\theta-3)/2}}.
\end{equation}}}
\begin{proof}{of Claim 2 :}
\begin{eqnarray*}
\E\left(\sum_{(j,k)\in\hat{m}}(P_n-P_n^*)^2(\psi_{j,k})\right)&\leq&\E\left(\sup_{m\in\M_n}\sum_{(j,k)\in m}(P_n-P_n^*)^2(\psi_{j,k})\right)\\
&\leq&\sum_{m\in\M_n}\sum_{(j,k)\in m}\E\left((P_n-P_n^*)^2(\psi_{j,k})\right)\\
&\leq&\frac{2}{p^2}\sum_{m\in\M_n}\sum_{l,l'=1}^{p}(g_{A,m}(j,k,l,l')+g_{B,m}(j,k,l,l'))
\end{eqnarray*}
with $$g_{m,A}(j,k,l,l')=\E\left( \sum_{(j,k)\in m}\left(\bar{\psi}_{j,k}(A_l)-\bar{\psi}_{j,k}(A_l^*)\right)\left(\bar{\psi}_{j,k}(A_{l'})-\bar{\psi}_{j,k}(A_{l'}^*)\right)\right) .$$
We develop this last term and we get, since
$$\left|\bar{\psi}_{j,k}(x)-\bar{\psi}_{j,k}(y)\right|\leq \frac{K_L2^{3j/2}\left|x-y\right|_q}{2q}$$
\begin{eqnarray*}
g_{A,m}(j,k,l,l')&\leq&\E\left(\sum_{(j,k)\in m}\left|\bar{\psi}_{j,k}(A_l)-\bar{\psi}_{j,k}(A_l^*)\right|\left|\bar{\psi}_{j,k}(A_{l'})-\bar{\psi}_{j,k}(A_{l'}^*)\right|\right)\\
&\leq&\E\left(\sum_{(j,k)\in m}\left|\bar{\psi}_{j,k}(A_l)-\bar{\psi}_{j,k}(A_l^*)\right|K_L2^{3j/2}\frac{\left|A_{l'}-A_{l'}^*\right|_q}{2q}\right)\\
&\leq&\frac{K_L\tau_q}{2}\sup_{x,y\in\R^q}\left\{\sum_{(j,k)\in m}2^{3j/2}\left|\bar{\psi}_{j,k}(x)-\bar{\psi}_{j,k}(y)\right|\right\}\\
&\leq&\frac{K_L\tau_q}{4}\sum_{j=0}^{J_m}2^{3j/2}\sup_{x,y\in\R}\left\{\sum_{k\in \Z}\left|\psi_{j,k}(x)-\psi_{j,k}(y)\right|\right\}\\
&\leq& \frac{2}{3}AK_LK_{\infty}2^{2J_m}\tau_q\;\textrm{since}\;\left\|\sum_{k\in\Z}|\psi_{j,k}|\right\|_{\infty}\leq AK_{\infty}2^{j/2}
\end{eqnarray*}
We can do the same computations for the term $g_{B,m}(j,k,l,l')$ and we obtain
$$\E\left(\sum_{j,k\in\hat{m}}\left((P_n-P_n^*)(\psi_{j,k})\right)^2\right)\leq L\tau_q\sum_{m\in \M_n}2^{2J_m}\leq L\tau_q2^{2J_n}\leq L\frac{(\log n)^{\kappa(\theta+1)}}{n^{(\theta-3)/2}}.$$
The last inequality comes from $q\geq \sqrt{n}/(2(\log n)^{\kappa})$ and Assumption {\bf{[AR]}}, the one before comes from Assumption {\bf{[W]}}. $\square$
\end{proof}
{\bf{Claim 3.}} {\it{Let us keep the notations of Theorem \ref{tauthm}, let $u=6/(7+\theta)<1/2$ and recall that $\kappa>2$. Let $\gamma$ be a real number in $(1,\kappa/2)$. Let
$$L_1^2=AK_{\infty}K_{BV}\sum_{l=0}^{\infty}\tilde{\beta}_l,\;L^2_2=2\Phi K^u_{BV}\sum_{k=0}^{\infty}\tilde{\beta}^{u}_k,\;L_3=\kappa(\epsilon)\Phi$$
\begin{equation}\label{L_1m}
\textrm{and}\;L_{1,m}=4(1+\epsilon)\left((1+\epsilon)L_1+L_2\sqrt{\frac{(\log D_m)^{\gamma}}{D_m^{1/2-u}}}+L_3\frac{(\log D_m)^{\gamma}}{(\log n)^{\kappa}}\right)^2,
\end{equation}
There exists a constant $L_{s}$ such that
$$\E\left(\sup_{m\in\M_n}\left\{\sum_{(j,k)\in m}(\nu_n^*)^2(\psi_{j,k})-\frac{L_{1,m}D_m}{n}\right\}\right)\leq \frac{L_{s}}{n}.$$}}
{\bf Remark :}
The series $\sum_{l=0}^{\infty}\tilde{\beta}_l$ and $\sum_{k=0}^{\infty}\tilde{\beta}^{u}_k$ are convergent under our hypotheses on the coefficients $\tau$. Since $s\in L^2([0,1])$, we have from Inequality (\ref{cvg}), $\tilde{\beta}_l\leq 2\|s\|_2^{2/3}\tau_l^{1/3}$ and thus $\tilde{\beta}_l\leq 2 \|s\|_2^{2/3}(1+l)^{-(1+\theta)/3}$. The series $\sum_{k=0}^{\infty}\tilde{\beta}^{u}_k$ converge since $\theta>5$ and 
$$\frac{u(1+\theta)}{3}=\frac{2(1+\theta)}{7+\theta}=1+\frac{\theta-5}{\theta+7}>1.$$ We use here $\tilde{\beta}$ instead of $\tau$ which allows to take $L_1$ not depending on $\|s\|_2$.

\begin{proof}{of Claim 3 :}
As in the previous section we use the following decomposition
\begin{eqnarray*}
\sum_{(j,k)\in m}(\nu_n^*)^2(\psi_{j,k})&=&\sum_{(j,k)\in m}\left(\bar{\nu}_{A,p}(\bar{\psi}_{j,k})+\bar{\nu}_{B,p}(\bar{\psi}_{j,k})\right)^2\\
&\leq&2\sum_{(j,k)\in m}\left(\bar{\nu}_{A,p}(\bar{\psi}_{j,k})\right)^2+2\sum_{(j,k)\in m}\left(\bar{\nu}_{B,p}(\bar{\psi}_{j,k})\right)^2
\end{eqnarray*}
We treat both terms with Proposition \ref{Massart} applied to the random variables $(A^*_l)_{0=1,..,p-1}$ and $(B^*_l)_{l=0,..,p-1}$ and to the class of functions $ \left\lbrace (\bar{\psi}_{j,k})_{(j,k)\in m}\right\rbrace $. Let
$$B_m^2=\sum_{(j,k)\in m}\Var\left(\bar{\psi}_{j,k}(A_1)\right),\;V_m^2= \sup_{t\in B_1(S_m)}\Var(\bar{t}(A_1)),\; H_m^2=\Vert\sum_{(j,k)\in m}\bar{\psi}_{j,k}^2\Vert_{\infty}.$$
We have, from Proposition \ref{Massart} 
\begin{equation}\label{conc}
\forall x>0,\;\p\left[\sqrt{\sum_{(j,k)\in m}(\bar{\nu}_{A,p})^2\bar{\psi}_{j,k}}\geq\frac{(1+\epsilon)}{\sqrt{p}}B_m+V_m\sqrt{\frac{2x}{p}}+\kappa(\epsilon)\frac{H_mx}{p}\right]\leq e^{-x}.
\end{equation}
Let us now evaluate $B_m$, $V_m$ and $H_m$, we have 
$$B_m^2=\frac{1}{(2q)^2}\sum_{(j,k)\in m}\Var\left(\sum_{i=1}^q\psi_{j,k}(X_i)\right).$$
From (\ref{bv}) and (\ref{inf}) we have $\forall j,k$ $\left\|\psi_{j,k}\right\|_{BV}\leq K_{BV}2^{j/2}$ and $\forall j$ $\Vert\sum_{k\in \Z}|\psi_{j,k}|\Vert_{\infty}\leq AK_{\infty}2^{j/2}$. Thus, from Inequality (\ref{covT})
\begin{eqnarray*}
\sum_{(j,k)\in m} \textrm{Var}\left(\sum_{i=1}^q\psi_{j,k}(X_i)\right)&\leq&2\sum_{(j,k)\in m} \sum_{l=1}^q(q+1-l)|\textrm{Cov}(\psi_{j,k}(X_1),\psi_{j,k}(X_l))|\\
&\leq& 2q\sum_{j=0}^{J_m}\sum_{k\in\Z}\sum_{l=1}^q\left\|\psi_{j,k}\right\|_{BV}\E\left(|\psi_{j,k}(X_1)|b(\sigma(X_1),X_l)\right)\\
&\leq& 2K_{BV}q\sum_{j=0}^{J_m}2^{j/2}\left\|\sum_{k\in \Z}|\psi_{j,k}(X_0)|\right\|_{\infty}\sum_{l=1}^q\tilde{\beta}_{l-1}\\
&\leq& 2q\left(AK_{\infty}K_{BV}\sum_{l=0}^{\infty}\tilde{\beta}_l\right)D_m.
\end{eqnarray*}
The last inequality comes from Assumption {\bf{[W]}}. \\
Since $L_1^2=AK_{\infty}K_{BV}\sum_{l=0}^{\infty}\tilde{\beta}_l$ we have 
\begin{equation}\label{BM}
B_m^2\leq\frac{ L_1^2D_m}{2q}.
\end{equation}
Let us deal with the term $V_m^2$. We have
\begin{equation}\label{Vm}
V_m^2\leq\sup_{t\in B_1(S_m)}\Var(\bar{t}(A_1))\leq\frac{2}{(2q)^2}\sum_{k=1}^{q}(q+1-k)\sup_{t\in B_1(S_m)}|\textrm{Cov}(t(X_1),t(X_k))|
\end{equation}
From Inequality (\ref{covT}), we have 
\begin{equation*}
|\Cov(t(X_1),t(X_k))|\leq \left\|t\right\|_{BV}\left\|t\right\|_{\infty}\tilde{\beta}_{k-1}.
\end{equation*}
Since $t$ belongs to $B_1(S_m)$, we have $t=\sum_{(j,k)\in m}a_{j,k}\psi_{j,k}$, with $\sum_{(j,k)\in m}a_{j,k}^2\leq 1$. Thus, by Cauchy-Schwarz inequality
\begin{eqnarray*}
\sum_{i=1}^l\left|t(x_{i+1})-t(x_i)\right|&\leq&\sum_{(j,k)\in m}|a_{j,k}|\sum_{i=1}^l\left|\psi_{j,k}(x_{i+1})-\psi_{j,k}(x_i)\right|\\
&\leq& \left(\sum_{(j,k)\in m}a_{j,k}^2\right)^{1/2}\left(\sum_{(j,k)\in m}\left(\sum_{i}|\psi_{j,k}(x_{i+1})-\psi_{j,k}(x_i)|\right)^2\right)^{1/2}\\
&\leq& \left(\sum_{(j,k)\in m}\left\|\psi_{j,k}\right\|^2_{BV}\right)^{1/2}\leq K_{BV}D_m. 
\end{eqnarray*}
Thus $\left\|t\right\|_{BV}\leq D_mK_{BV}$. From Assumption {\bf{[$M_2$]}}, we have $\left\|t\right\|_{\infty}\leq \Phi\sqrt{D_m}$. Thus 
\begin{equation}\label{co1}
|\textrm{Cov}(t(X_1),t(X_k))|\leq \Phi K_{BV}\tilde{\beta}_{k-1}D_m^{3/2}.
\end{equation}
Moreover, we have by Cauchy-Schwarz inequality and {\bf{[$M_2$]}}
\begin{equation}\label{co2}
|\textrm{Cov}(t(X_1),t(X_k))|\leq \left\|t\right\|_{\infty}\left\|t\right\|_2\left\|s\right\|_2\leq \Phi\left\|s\right\|_2\sqrt{D_m}.
\end{equation}
We use the inequality $a\wedge b\leq a^{u}b^{1-u}$ with 
$$a=\Phi K_{BV}\tilde{\beta}_{k-1}D_m^{3/2},\;b=\Phi\left\|s\right\|_2\sqrt{D_m},\;u=\frac6{7+\theta}<\frac12.$$
From (\ref{co1}) and (\ref{co2}), we derive that
$$|\textrm{Cov}(t(X_1),t(X_k))|\leq L_k'D_m^{1/2+u}\;\textrm{where} \;L_k'=\Phi\left(K_{BV}\tilde{\beta}_{k-1}\right)^{u}\left\|s\right\|^{1-u}_2.$$
Pluging this inequality in (\ref{Vm}), we obtain
\begin{equation}\label{VM}
V_m^2\leq \frac{L_2^2\left\|s\right\|^{1-u}_2D_m^{1/2+u}}{4q}\;\textrm{since} \;L_2^2=2\Phi K_{BV}^{u}\sum_{k=0}^{\infty}\tilde{\beta}^{u}_k.
\end{equation}

Finally, we have from hypothesis {\bf{[$M_2$]}} 
\begin{equation}\label{HM}
H_m^2\leq \frac{1}{4}\left\|\sum_{(j,k)\in m}\psi_{j,k}^2\right\|_{\infty}\leq \frac{\Phi^2 D_m}{4}.
\end{equation}
Let $y>0$ and let us apply Inequality (\ref{conc}) with $x=((\log D_m)^{\gamma}/\left\|s\right\|^{1-u}_2)+(y/D_m^{1/2+u})$. We have, from (\ref{BM}), (\ref{VM}) and (\ref{HM})
\begin{eqnarray*}
&&\p\left[\sum_{(j,k)\in m}(\bar{\nu}_{A,p})^2(\bar{\psi}_{j,k})
>\left((1+\epsilon)\sqrt{\frac{L_1^2D_m}{2pq}}+\frac{L_3\sqrt{D_m}}{2p}\left(\frac{(\log D_m)^{\gamma}}{\|s\|^{1-u}_2}+\frac{y}{D_m^{1/2+u}}\right)\right.\right.\\
&&\left.\left.+\sqrt{\frac{L_2^2\|s\|^{1-u}_2D_m^{1/2+u}}{2pq}\left(\frac{(\log D_m)^{\gamma}}{\|s\|^{1-u}_2}+\frac{y}{D_m^{1/2+u}}\right)}\right)^2 \right]\leq e^{-\frac{(\log D_m)^{\gamma}}{\|s\|^{1-u}_2}}e^{-D_m^{-(1/2+u)}y}.
\end{eqnarray*}
Then, we use the inequality $\sqrt{\alpha+\beta}\leq \sqrt{\alpha}+\sqrt{\beta}$ with 
$$\alpha= \frac{(\log D_m)^{\gamma}}{\|s\|^{1-u}_2}\;\textrm{and}\; \beta=\frac{y}{D_m^{1/2+u}}$$
and the inequality $(a+b)^2\leq (1+\epsilon)a^2+(1+\epsilon^{-1})b^2$ with  
$$a=\left((1+\epsilon)L_1+L_2\sqrt{\frac{(\log D_m)^{\gamma}}{D_m^{1/2-u}}}
+\frac{L_3(\log D_m)^{\gamma}}{\left\|s\right\|^{1-u}_2(\log n)^{\kappa}}\right)\sqrt{\frac{D_m}{n}}$$
$$\textrm{and}\;b=\frac{1}{\sqrt{n}}\left(L_2\sqrt{\|s\|^{1-u}_2y}+\frac{L_3 y}{(\log n)^{\kappa}D_m^{u}}\right).$$ 
Setting $L_m=(1+\epsilon)a^2n/D_m$, we obtain
\begin{eqnarray*}
&&\p\left(\sum_{(j,k)\in m}(\bar{\nu}_{A,p})^2(\bar{\psi}_{j,k})-\frac{L_mD_m}{n}>\frac{(1+\epsilon^{-1})}{n}\left(L_2\sqrt{\|s\|^{1-u}_2y}+\frac{L_3 y}{(\log n)^{\kappa}D_m^{u}}\right)^2\right)\\
&&\leq e^{-\frac{(\log D_m)^{\gamma}}{\|s\|^{1-u}_2}}e^{-D_m^{-(1/2+u)}y}.
\end{eqnarray*}
Thus, for all $y>0$,
$$\p\left(\sup_{m\in\M_n}\left\{\sum_{(j,k)\in m}(\bar{\nu}_{A,p})^2(\bar{\psi}_{j,k})- \frac{L_mD_m}{n}\right\}>\frac{L_s}{n}(y+y^2)\right)\leq \sum_{m\in\M_n}e^{-\frac{(\log D_m)^{\gamma}}{\left\|s\right\|^{1-u}_2}-D_m^{-(1/2+u)}y}$$
where $L_s=2(1+\epsilon^{-1})\left[(L_2\sqrt{\|s\|^{1-u}_2})\vee L_3 /((\log 2)^{\kappa}2^{u})\right]^2$.
We can integrate this last inequality to prove Claim 3.$\square$
\end{proof}
{\bf{Claim 4 :}}{\it{We keep the notations of the previous Claims. Let
\begin{equation}\label{L2mm}
L_2(m,m')=4\left(L_2\sqrt{\frac{(\log (D_m\vee D_{m'}))^{\gamma}}{(D_m\vee D_{m'})^{1/2-u}}}+\frac{\Phi}{3(\log n)^{\kappa-\gamma}}\right)^2.
\end{equation}
Then there exists a constant $L_{s, \theta}$ depending on $\left\|s\right\|_2$  and $\theta$ such that, for all $\eta>0$ 
$$\E\left(\sup_{m,m'\in\M_n}\left\{ \nu_n(s_{m}-s_{m'})-\frac{\left\|s_{m}-s_{m'}\right\|_2^2}{2\eta}-\eta\frac{L_2(m,m')(D_{m}\vee D_{m'})}{n}\right\}\right)\leq \frac{\eta L_{s,\theta}}{n}.$$}}
\begin{proof}{of Claim 4 :}
\begin{eqnarray}\label{decB}
&&\E\left(\sup_{m,m'\in\M_n}\left\{ \nu_n(s_{m}-s_{m'})-\frac{\left\|s_{m}-s_{m'}\right\|_2^2}{2\eta}-\eta\frac{L_2(m,m')(D_{m}\vee D_{m'})}{n}\right\}\right)\nonumber\\
&&\leq \E\left(\sup_{m,m'}(P_n-P_n^*)(s_{m}-s_{m'})\right)\nonumber\\
&&+\E\left(\sup_{m,m'}\left \{\nu_n^*(s_{m}-s_{m'})-\frac{\left\|s_{m}-s_{m'}\right\|_2^2}{2\eta}-\eta\frac{L_2(m,m')(D_{m}\vee D_{m'})}{n}\right\}\right).
\end{eqnarray}
Since $\forall l=0,...,p-1$, $\E\left(|A_l-A_l^*|_q\right)\leq q\tau_q$, we have 
\begin{eqnarray}
\E\left(\sup_{m,m'}(P_n-P_n^*)(s_{m}-s_{m'})\right)&\leq& 2\sum_{m,m'}\E \left(|(\bar{s}_m-\bar{s}_{m'})(A_1)-(\bar{s}_m-\bar{s}_{m'})(A_1^*)|\right)\nonumber\\
&\leq& \tau_q\sum_{m,m'}\textrm{Lip}(s_m-s_{m'}).\nonumber
\end{eqnarray}
When $m\subset m'$, we have, for all $x,y\in \R$, using Assumption {\bf{[W]}},
\begin{equation*}
\frac{|(s_m-s_{m'})(x-y)|}{|x-y|}\leq \sum_{j=J_m+1}^{J_{m'}}\sum_{k=-A_2}^{2^{j}-A_1}|P\psi_{j,k}| \frac{|\psi_{j,k}(x)-\psi_{j,k}(y)|}{|x-y|}
\end{equation*}
Let us fix $j\in [J_m+1,J_{m'}]$, from Assumption {\bf{[W]}}, there is less than $A$ indexes $k\in \Z$ such that $\psi_{j,k}(x)\neq 0$, thus there is less than $2A$ indexes such that $|\psi_{j,k}(x)-\psi_{j,k}(y)|\neq 0$. Hence 
\begin{eqnarray*}
\sum_{k\in \Z}|P\psi_{j,k}|\frac{|\psi_{j,k}(x)-\psi_{j,k}(y)|}{|x-y|}&\leq& 2A\sup_{k\in\Z}|P\psi_{j,k}|\Lip(\psi_{j,k})\\
&\leq&2A\left\|s\right\|_2K_L2^{3j/2}. 
\end{eqnarray*}
Thus, $\textrm{Lip}(s_m-s_{m'})\leq A\left\|s\right\|_2K_L\sqrt{8}2^{3J_{m'}/2}/(\sqrt{8}-1)$ and by Assumptions {\bf{[W]}}, {\bf{[AR]}} and the value of $q$, 
\begin{equation}\label{rb}
\E\left(\sup_{m,m'}(P_n-P_n^*)(s_{m}-s_{m'})\right)\leq L_sn^{3/2}(\log n)\tau_q\leq L_s\frac{(\log n)^{\kappa(\theta+1)+1}}{n^{(\theta-2)/2}}.
\end{equation}
Let us deal with the other term in (\ref{decB}). We have, $\forall \eta>0$
\begin{eqnarray}
\nu_n^*(s_{m}-s_{m'})&\leq& \frac{\Vert s_{m}-s_{m'}\Vert_2^2}{2\eta}+\frac{\eta}{2}\left(\bar{\nu}_{A,p}(\bar{t}_{m,m'})+\bar{\nu}_{B,p}(\bar{t}_{m,m'})\right)^2 \nonumber \\
&\leq& \frac{\Vert s_{m}-s_{m'}\Vert_2^2}{2\eta}+ \eta ( \bar{\nu}_{A,p}(\bar{t}_{m,m'}))^2+\eta ( \bar{\nu}_{B,p}(\bar{t}_{m,m'}))^2 \label{bern}
\end{eqnarray}
where, as in the proof of Theorem \ref{betathm}, $t_{m,m'}=(s_{m}-s_{m'})/\Vert s_{m}-s_{m'}\Vert_2$. We apply Bernstein's inequality to the function $\bar{t}_{m,m'}$ and the variables $A_l^*$, we have
\begin{equation}\label{berni}
\forall x>0,\;\p\left(\bar{\nu}_{A,p}(\bar{t}_{m,m'})>\sqrt{\frac{2\Var(\bar{t}_{m,m'}(A_0))x}{p}}+\frac{\Vert \bar{t}_{m,m'}\Vert_{\infty}x}{3p}\right)\leq e^{-x}.
\end{equation}
We proceed as in the proof of Claim 3 to control this variance. We have, by stationarity of the process $(X_n)_{n\in\Z}$, 
$$\Var(\bar{t}_{m,m'}(A_0))=\frac{1}{2q^2}\sum_{k=0}^{q-1}(q-k)\textrm{Cov}(t_{m,m'}(X_1),t_{m,m'}(X_{k+1})).$$
From Inequality (\ref{covT}), we have
\begin{equation*}
\left|\Cov(t_{m,m'}(X_1),t_{m,m'}(X_{k+1}))\right|\leq \left\| t_{m,m'}\right\|_{BV}\left\| t_{m,m'}\right\|_{\infty} \tilde{\beta}_k.
\end{equation*}
Let $m\bigtriangleup m'$ be the set of indexes that belong to $m\cup m'$ but do not belong to $m\cap m'$. We use the same computations as in the proof of Claim 3 to get 
\begin{equation*}
\left\| t_{m,m'}\right\|_{BV}\leq \frac{\left\|\sum_{(j,k)\in m'\bigtriangleup m}(P\psi_{j,k})\psi_{j,k}\right\|_{BV}}{\left\|s_m-s_{m'}\right\|_2}\leq \sqrt{\sum_{(j,k)\in m'\bigtriangleup m}\left\|\psi_{j,k}\right\|_{BV}^2}\leq K_{BV}(D_m\vee D_{m'}).
\end{equation*}
Since $\left\|t_{m,m'}\right\|_{\infty}=\Phi\sqrt{D_m\vee D_{m'}}$, we have 
\begin{equation}\label{coB1}
\left|\textrm{Cov}(t_{m,m'}(X_1),t_{m,m'}(X_{k+1}))\right|\leq \Phi K_{BV}\tilde{\beta}_k(D_m\vee D_{m'})^{3/2}.
\end{equation}
Moreover, we have 
\begin{equation}\label{coB2}
\textrm{Cov}(t_{m,m'}(X_1),t_{m,m'}(X_{k+1}))\leq \left\|t_{m,m'}\right\|_{\infty}\left\|t_{m,m'}\right\|_2\left\|s\right\|_2\leq \Phi\left\|s\right\|_2\sqrt{(D_m\vee D_m')}.
\end{equation}
Thus, using $a\wedge b\leq a^{u}b^{1-u}$ with 
$$a=\Phi K_{BV}\tilde{\beta}_k(D_m\vee D_{m'})^{3/2},\;b=\Phi\left\|s\right\|_2\sqrt{(D_m\vee D_{m'})},\;\textrm{and}\;u=\frac6{7+\theta}<\frac12,$$
we have
\begin{equation*}
\left|\textrm{Cov}(t_{m,m'}(X_1),t_{m,m'}(X_{k+1}))\right|\leq\Phi K_{BV}^u\tilde{\beta}^{u}_k\left\|s\right\|^{1-u}_2(D_m\vee D_{m'})^{1/2+u}. 
\end{equation*}
Thus 
\begin{equation}\label{variance}
\Var(\bar{t}_{m,m'}(A_0))\leq\Phi K_{BV}^u\left(\sum_{k=0}^{\infty}\tilde{\beta}^{u}_k\right)\left\|s\right\|^{1-u}_2\frac{(D_m\vee D_{m'})^{1/2+u}}{2q}.
\end{equation}
Moreover 
\begin{equation}\label{supnorm}
\Vert \bar{t}_{m,m'}\Vert_{\infty}\leq \frac12\Vert t_{m,m'}\Vert_{\infty}\leq \frac12\Phi\sqrt{D_m\vee D_m'}.
\end{equation}
Now, we use (\ref{berni}) with $x=(\log (D_m\vee D_{m'}))^{\gamma}/\left\|s\right\|^{1-u}_2+y/(D_m\vee D_{m'})^{1/2+u}$. From (\ref{variance}) and (\ref{supnorm}), we have for all $y>0$,
\begin{eqnarray*}
&&\p\left(\bar{\nu}_{A,p}(\bar{t}_{m,m'})>L_2\sqrt{\frac{(D_m\vee D_{m'})^{1/2+u}}{2pq}\left((\log (D_m\vee D_{m'}))^{\gamma}+\frac{\left\|s\right\|^{1-u}_2y}{(D_m\vee D_{m'})^{1/2+u}}\right)}\right.\\
&&\left.+\frac{\Phi\sqrt{D_m\vee D_m'}}{6p}\left(\frac{(\log (D_m\vee D_{m'}))^{\gamma}}{\left\|s\right\|^{1-u}_2}+\frac{y}{(D_m\vee D_{m'})^{1/2+u}}\right)\right)\\
&&\leq e^{-\frac{(\log (D_m\vee D_{m'}))^{\gamma}}{\left\|s\right\|^{1-u}_2}}e^{-\frac{y}{(D_m\vee D_{m'})^{1/2+u}}}.
\end{eqnarray*}
Now we use the inequality $\sqrt{a+b}\leq \sqrt{a}+\sqrt{b}$ with 
$$a=(\log (D_m\vee D_{m'}))^{\gamma}\;\textrm{and}\;b=\frac{\left\|s\right\|^{1-u}_2y}{(D_m\vee D_{m'})^{1/2+u}}$$
and we obtain, using Assumption {\bf{[$M_1$]}} 
\begin{eqnarray*}
&&\p\left(\bar{\nu}_{A,p}\bar{t}_{m,m'}-\sqrt{\frac{L_2(m,m')(D_m\vee D_m')}{n}}>\frac {L_s}{\sqrt{n}}(\sqrt{y}+y)\right)\\
&&\leq e^{-\frac{(\log (D_m\vee D_{m'}))^{\gamma}}{\left\|s\right\|^{1-u}_2}}e^{-(D_m\vee D_{m'})^{-(1/2+u)}y},
\end{eqnarray*}
with 
$$L_2(m,m')=\left(L_2\sqrt{\frac{(\log (D_m\vee D_{m'}))^{\gamma}}{(D_m\vee D_{m'})^{1/2-u}}}+\frac{\Phi(\log(D_m\vee D_{m'}))^{\gamma}}{3(\log n)^{\kappa}}\right)^2,$$
$$\textrm{and}\;L_s=L_2\sqrt{\|s\|^{1-u}_2}\vee \frac{\Phi}{3(\log 2)^{\kappa}2^{u}}.$$
Thus, we obain
\begin{eqnarray*}
&&\p\left((\bar{\nu}_{A,p}\bar{t}_{m,m'})^2>2\frac{L_2(m,m')(D_m\vee D_m')}{n}+4\frac{L_s^2}{n}(y+y^2)\right) \\
&&\leq e^{-\frac{(\log (D_m\vee D_{m'}))^{\gamma}}{\left\|s\right\|^{1-u}_2}-\frac{y}{(D_m\vee D_{m'})^{1/2+u}}}.
\end{eqnarray*}
The same result holds for $\bar{\nu}_{B,p}\bar{t}_{m,m'}$. Thus we obtain from (\ref{bern})
\begin{eqnarray*}
&&\p\left(\nu_n^*(s_{m}-s_{m'})\geq \frac{\Vert s_{m}-s_{m'}\Vert_2^2}{2\eta}+4\eta\frac{L_2(m,m')(D_m\vee D_m')}{n}+8\eta\frac{L_s^2}{n}(y+y^2)\right)\\
&&\leq 2e^{-\frac{(\log (D_m\vee D_{m'}))^{\gamma}}{\left\|s\right\|^{1-u}_2}-\frac{y}{(D_m\vee D_{m'})^{1/2+u}}}.
\end{eqnarray*}
We deduce that 
\begin{eqnarray*}
&&\p\left(\exists m,m'\in\M_n,\;\nu_n^*(s_{m}-s_{m'})- \frac{\Vert s_{m}-s_{m'}\Vert_2^2}{2\eta}-4\eta\frac{L_2(m,m')(D_m\vee D_m')}{n}\right.\\
&&\left.\geq8\eta\frac{L_s^2}{n}(y+y^2)\right)\leq 2\sum_{m,m'\in\M_n}\left(e^{-\frac{(\log (D_m\vee D_{m'}))^{\gamma}}{\left\|s\right\|^{1-u}_2}}\right)e^{-\frac{y}{(D_m\vee D_{m'})^{1/2+u}}}.
\end{eqnarray*}
We integrate this last inequality to get Claim 4.$\square$ 
\end{proof}
{\bf{Conclusion of the proof:}}\\
Take
$$\pen'(m)\geq (2L_{1,m}+\eta L_2(m,m))\frac{D_m}{n},$$ 
where $L_{1,m}$ and $L_2(m,m)$ are defined by (\ref{L_1m}) and (\ref{L2mm}) respectively.
From Claims 2, 3 and 4, if we take the expectation in (\ref{dec}), we have, for some constant $L_s$,
\begin{equation}\label{tp}
\E\left(\left\|s-\tilde{s}\right\|_2^2\right)\leq \E\left(\left\|s-\hat{s}_{m_o}\right\|_2^2+\pen'(m_o)-V(m_o)+2\eta L_2(m_o,m_o)\frac{D_{m_o}}{n}\right)+\frac{\eta L_s}{n}.
\end{equation}
Moreover, if $D_m\geq \left((L_2/L_1)(\log n)^{\kappa-\gamma/2}\right)^{2(7+\theta)/(\theta-5)}$, we have
\begin{eqnarray}
\frac{L_{1,m}}{4L_1^2 }&\leq& (1+\epsilon)\left((1+\epsilon)+\left(1+\frac{L_3}{2L_1}\right)(\log n)^{-(\kappa-\gamma)}\right)^2\nonumber\\
&\leq& (1+\epsilon)^3+(1+\epsilon^{-1})(1+\epsilon)\left(1+\frac{L_3}{2L_1}\right)^2(\log n)^{-2(\kappa-\gamma)}\label{step}.
\end{eqnarray}
We use the inequality $(a+b)^2\leq (1+\epsilon)a^2+(1+\epsilon^{-1})b^2$ to obtain (\ref{step}). Moreover, we have 
$$L_2(m,m)\leq 4L_1^2\left(\left(1+\frac{\Phi}{6L_1}\right)(\log n)^{-(\kappa-\gamma)}\right)^2.$$
As in the proof of Theorem \ref{betathm}, we take $\eta=(\log n)^{\kappa-\gamma}$ and we fix $\epsilon$ sufficiently small. For $n\geq n_o$, we have $2L_{1,m}+\eta L_2(m,m)<KL_1^2$. Thus inequality (\ref{tau}) follows from (\ref{tp}).$\square$
\section{Appendix}
This section is devoted to technical lemmas that are needed in the proofs. 
\subsection{Covariance inequality}
\begin{lemma}{Viennet's inequality}
Let $(X_n)_{n\in \Z}$ be a stationary and $\beta$-mixing process. There exists a positive function $b$ such that $P(b)\leq\sum_{l=0}^{\infty}\beta_l$, $P(b^p)\leq p\sum_{l=1}^{\infty}l^{p-1}\beta_l$, and for all function $h\in L_2(P)$ 
\begin{equation}\label{Viennet}
\textrm{Var}\left(\sum_{l=1}^q h(X_l)\right)\leq 4qP(bh^2).
\end{equation}
\end{lemma}
\subsection{Concentration inequalities}
We sum up in this section the concentration inequalities we used in the proofs. We begin with Bernstein's inequality
\begin{prop}{Bernstein's inequality}\\
Let $X_1,...,X_n$ be iid random variables valued in a measurable space $(X,\mathcal{X})$ and let $t$ be a measurable real valued function. Let $v=\Var (t(X_1))$ and $b=\left\|t\right\|_{\infty}$, then, for all $x>0$, we have
$$\p\left((P_n-P)t>v\sqrt{\frac{2x}{n}}+\frac{bx}{3n}\right)\leq e^{-x}.$$
\end{prop}
Now we give the most important tool of our proof, it is a concentration's inequality for the supremum of the empirical process over a class of function. We give here the version of Bousquet \cite{Bo02}.
\begin{theo}{Talagrand's Theorem}\\
Let $X_1,...,X_n$ be i.i.d random variables valued in some measurable space $[X,\mathcal{X}]$. Let $\F$ be a separable class of bounded functions from $X$ to $\R$ and assume that all functions $t$ in $\F$ are $P$-measurable, and satisfy $\Var(t(X_1))\leq \sigma^2$, $\Vert t\Vert_{\infty}\leq b$. Then 
$$\p\left(\sup_{t\in \F}\nu_n(t)>\E\left(\sup_{t\in \F}\nu_n(t)\right)+\sqrt{\frac{2x\left(\sigma^2+2b\E\left(\sup_{t\in \F}\nu_n(t)\right)\right)}{n}}+\frac{bx}{3n}\right)\leq e^{-x}.$$
In particular, for all $\epsilon>0$, if $\kappa(\epsilon)=1/3+\epsilon^{-1}$, we have
$$\p\left(\sup_{t\in \F}\nu_n(t)>(1+\epsilon)\E\left(\sup_{t\in \F}\nu_n(t)\right)+\sigma\sqrt{\frac{2x}{n}}+\kappa(\epsilon)\frac{bx}{n}\right)\leq e^{-x}.$$
\end{theo}
We can deduce from this Theorem a concentration's inequality for $\chi$-square type statistics. This is Proposition (7.3) of Massart \cite{Ma07}.
\begin{prop}\label{Massart}
Let $X_1,...,X_n$ be independent and identically distributed random variables valued in some measurable space $(X,\mathcal{X})$. Let $P$ denote their common distribution. Let $\phi_{\lambda}$ be a finite family of measurable and bounded functions on $(X,\mathcal{X})$. Let
$$H_{\Lambda}^2=\Vert\sum_{\lambda\in\Lambda}\phi_{\lambda}^2\Vert_{\infty}\;\textrm{and}\;B^2_{\Lambda}=\sum_{\lambda\in\Lambda}\Var(\phi_{\lambda}(X_1)).$$
Moreover, let $\mathcal{S}_{\Lambda}=\left\{a\in \R^{\Lambda}:\sum_{\lambda\in\Lambda}a_{\lambda}^2= 1\right\}$ and 
$$V_{\Lambda}^2=\sup_{a\in\mathcal{S}_{\Lambda}}\left\{\Var\left(\sum_{\lambda\in\Lambda}a_{\lambda}\phi_{\lambda}(X_1)\right)\right\}.$$
Then the following inequality holds, for all positive $x$ and $\epsilon$
\begin{equation}\label{massart}
\p\left[\left(\sum_{\lambda\in \Lambda}(P_n-P)^2\phi_{\lambda}\right)^{1/2}\geq \frac{1+\epsilon}{\sqrt{n}}B_{\Lambda}+V_{\Lambda}\sqrt{\frac{2x}{n}}+\kappa(\epsilon)\frac{H_{\Lambda}x}{n}\right]\leq e^{-x},
\end{equation}
where $\kappa(\epsilon)=\epsilon^{-1}+1/3$.
%Let $\epsilon>0$ be given. Consider the set
%$$C_{\Lambda}=\left\lbrace a; \sup_{\lambda}|a_{\lambda}|=1\right\rbrace $$
%and define
%$$\sigma_P^2=\sup_{\sum a_{\lambda}^2\leq 1}\left[ V_P(\sum a_{\lambda}\phi_{\lambda})\right]\;\textrm{and}\;H=\sup_{a\in C_{\Lambda}}\left\|\sum a_{\lambda}\phi_{\lambda}\right\|_{\infty}. $$
%Moreover, for any subset $m$ of $\Lambda$ let 
%$$B_m^2:=\frac{\sum_{\lambda\in m}V_P(\phi_{\lambda})}{n}.$$
%Then setting $\kappa(\epsilon)=2(\epsilon^{-1}+1/3)$, $\tilde{\kappa}(\epsilon)=2(\epsilon^{-1}+8/3)$, there exists events $\Omega_n$ and $\tilde{\Omega}_n$ such that on one hand
%$$\p(\Omega_n^c)\leq 2|\Lambda|\exp\left(-\eta(\epsilon)\frac{n\sigma_p^2}{H^2}\right),\;\p(\tilde{\Omega}_n^c)\leq 2|\Lambda|\exp\left(-\tilde{\eta}(\epsilon)\frac{n\sigma_p^2}{H^2}\right)$$
%where $$\eta(\epsilon)=\frac{2\epsilon^2}{\kappa(\epsilon)(\kappa(\epsilon)+\frac{2\epsilon}{3})}$$
%and $\tilde{\eta}$ defined as $\eta$ but with $\tilde{\kappa}$in stead of $\kappa$ and, on the other hand, for any subset $m$ of $\gamma$ end any positive $x$
%$$\p\left[\left(\sum_{\lambda\in m}(P_n-P)^2\phi_{\lambda}\right)^{1/2}{\bf{1}}_{\Omega_n}\geq (1+\epsilon)\left(B_m+\sigma_P\sqrt{\frac{2x}{n}}\right)\right]\leq e^{-x}$$
%$$\p\left[\left(\sum_{\lambda\in m}(P_n-P)^2\phi_{\lambda}\right)^{1/2}{\bf{1}}_{\tilde{\Omega}_n}\leq (1+\epsilon)\left(B_m-\sigma_P\sqrt{\frac{2x}{n}}\right)\right]\leq e^{-x}.$$
\end{prop}
\begin{proof}:
Following Massart \cite{Ma07} Proposition 7.3, we remark that, by Cauchy-Schwarz's inequality
$$\left(\sum_{\lambda\in \Lambda}\nu_n^2\phi_{\lambda}\right)^{1/2}=\sup_{a\in \mathcal{S}_{\Lambda}}\sum_{\lambda\in \Lambda}a_{\lambda}\nu_n\phi_{\lambda}=\sup_{a\in \mathcal{S}_{\Lambda}}\nu_n\left(\sum_{\lambda\in \Lambda}a_{\lambda}\phi_{\lambda}\right).$$
Thus the result follows by applying Talagrand's Theorem to the class of functions
$$\F=\left\{t=\sum_{\lambda\in \Lambda}a_{\lambda}\phi_{\lambda};\;a\in\mathcal{S}_{\Lambda} \right\}.$$
\end{proof}
\bibliographystyle{plain}
\bibliography{bibliolerasle}

\begin{thebibliography}{10}

\bibitem{Ak73}
H.~Akaike.
\newblock Information theory and an extension of the maximum likelihood
  principle.
\newblock In {\em Second {I}nternational {S}ymposium on {I}nformation {T}heory
  ({T}sahkadsor, 1971)}, pages 267--281. Akad{\'e}miai Kiad{\'o}, Budapest,
  1973.

\bibitem{Ak70}
Hirotugu Akaike.
\newblock Statistical predictor identification.
\newblock {\em Ann. Inst. Statist. Math.}, 22:203--217, 1970.

\bibitem{An84}
Donald W.~K. Andrews.
\newblock Nonstrong mixing autoregressive processes.
\newblock {\em J. Appl. Probab.}, 21(4):930--934, 1984.

\bibitem{AM08}
S.~Arlot and P.~Massart.
\newblock Data-driven calibration of penalties for least squares regression.
\newblock {\em Submitted to Journal of Machine learning research}, 2008.

\bibitem{Ar08}
Sylvain Arlot.
\newblock Model selection by resampling penalization.
\newblock {\em hal-00262478}, 2008.

\bibitem{BCV01}
Y.~Baraud, F.~Comte, and G.~Viennet.
\newblock Adaptive estimation in autoregression or {$\beta$}-mixing regression
  via model selection.
\newblock {\em Ann. Statist.}, 29(3):839--875, 2001.

\bibitem{Be79}
Henry C.~P. Berbee.
\newblock {\em Random walks with stationary increments and renewal theory},
  volume 112 of {\em Mathematical Centre Tracts}.
\newblock Mathematisch Centrum, Amsterdam, 1979.

\bibitem{BM97}
Lucien Birg{\'e} and Pascal Massart.
\newblock From model selection to adaptive estimation.
\newblock In {\em Festschrift for {L}ucien {L}e {C}am}, pages 55--87. Springer,
  New York, 1997.

\bibitem{BM07}
Lucien Birg{\'e} and Pascal Massart.
\newblock Minimal penalties for {G}aussian model selection.
\newblock {\em Probab. Theory Related Fields}, 138(1-2):33--73, 2007.

\bibitem{Bo02}
Olivier Bousquet.
\newblock A {B}ennett concentration inequality and its application to suprema
  of empirical processes.
\newblock {\em C. R. Math. Acad. Sci. Paris}, 334(6):495--500, 2002.

\bibitem{Br07}
Richard~C. Bradley.
\newblock {\em Introduction to strong mixing conditions. {V}ol. 1}.
\newblock Kendrick Press, Heber City, UT, 2007.

\bibitem{CDT08}
F.~Comte, J.~Dedecker, and M.~L. Taupin.
\newblock Adaptive density deconvolution with dependent inputs.
\newblock {\em Math. Methods Statist.}, 17(2):87--112, 2008.

\bibitem{CM02}
Fabienne Comte and Florence Merlev{\`e}de.
\newblock Adaptive estimation of the stationary density of discrete and
  continuous time mixing processes.
\newblock {\em ESAIM Probab. Statist.}, 6:211--238 (electronic), 2002.
\newblock New directions in time series analysis (Luminy, 2001).

\bibitem{DDLLLP}
J\'er\^ome Dedecker, Paul Doukhan, Gabriel Lang, Jos{\'e}~Rafael {Le{\'o}n R.},
  Sana Louhichi, and Cl{\'e}mentine Prieur.
\newblock {\em Weak dependence: with examples and applications}, volume 190 of
  {\em Lecture Notes in Statistics}.
\newblock Springer, New York, 2007.

\bibitem{DP05}
J\'er\^ome Dedecker and Cl{\'e}mentine Prieur.
\newblock New dependence coefficients. {E}xamples and applications to
  statistics.
\newblock {\em Probab. Theory Related Fields}, 132(2):203--236, 2005.

\bibitem{DJKP}
David~L. Donoho, Iain~M. Johnstone, G{\'e}rard Kerkyacharian, and Dominique
  Picard.
\newblock Density estimation by wavelet thresholding.
\newblock {\em Ann. Statist.}, 24(2):508--539, 1996.

\bibitem{Do94}
Paul Doukhan.
\newblock {\em Mixing}, volume~85 of {\em Lecture Notes in Statistics}.
\newblock Springer-Verlag, New York, 1994.
\newblock Properties and examples.

\bibitem{GW09}
Ir\`ene Gannaz and Olivier Wintenberger.
\newblock Adaptive density estimation under dependence.
\newblock {\em forthcoming in ESAIM, Probab. and Statist.}, 2008.

\bibitem{Ma73}
C.L. Mallows.
\newblock Some comments on $c_p$.
\newblock {\em Technometrics}, 15:661--675, 1973.

\bibitem{Ma07}
Pascal Massart.
\newblock {\em Concentration inequalities and model selection}, volume 1896 of
  {\em Lecture Notes in Mathematics}.
\newblock Springer, Berlin, 2007.
\newblock Lectures from the 33rd Summer School on Probability Theory held in
  Saint-Flour, July 6--23, 2003, With a foreword by Jean Picard.

\bibitem{Pr07}
C.~Prieur.
\newblock Change point estimation by local linear smoothing under a weak
  dependence condition.
\newblock {\em Math. Methods Statist.}, 16(1):25--41, 2007.

\bibitem{Ri00}
Emmanuel Rio.
\newblock {\em Th{\'e}orie asymptotique des processus al{\'e}atoires faiblement
  d{\'e}pendants}, volume~31 of {\em Math{\'e}matiques \& Applications (Berlin)
  [Mathematics \& Applications]}.
\newblock Springer-Verlag, Berlin, 2000.

\bibitem{Ru82}
Mats Rudemo.
\newblock Empirical choice of histograms and kernel density estimators.
\newblock {\em Scand. J. Statist.}, 9(2):65--78, 1982.

\bibitem{Ta96}
Michel Talagrand.
\newblock New concentration inequalities in product spaces.
\newblock {\em Invent. Math.}, 126(3):505--563, 1996.

\bibitem{Vi97}
Gabrielle Viennet.
\newblock Inequalities for absolutely regular sequences: application to density
  estimation.
\newblock {\em Probab. Theory Related Fields}, 107(4):467--492, 1997.

\bibitem{RV59}
V.~A. Volkonski{\u\i} and Yu.~A. Rozanov.
\newblock Some limit theorems for random functions. {I}.
\newblock {\em Teor. Veroyatnost. i Primenen}, 4:186--207, 1959.

\end{thebibliography}

\end{document}